\documentclass[11pt,a4paper]{amsart}
\usepackage{amssymb,amsmath,amsthm}
\theoremstyle{plain}
\newtheorem{theorem}{Theorem}[section]
\newtheorem{definition}[theorem]{Definition}
\newtheorem{lemma}[theorem]{Lemma}
\newtheorem{corollary}[theorem]{Corollary}
\newtheorem{proposition}[theorem]{Proposition}

\theoremstyle{remark}
\newtheorem{remark}[theorem]{Remark}

\newtheorem{example}[theorem]{Example}

\usepackage[dvipdfm,
  bookmarks=true,
  bookmarksnumbered=false,
  bookmarkstype=toc]{hyperref}

\numberwithin{equation}{section}

\newcommand{\R}{\mathbb{R}}

\renewcommand{\Im}{\operatorname{Im}}

\newcommand{\I}{\infty}

\def\({\left(}
\def\){\right)}
\def\<{\left\langle}
\def\>{\right\rangle}
\def\le{\leqslant}
\def\ge{\geqslant}

\def\l{\lambda}

\newcommand{\Ga}{\Gamma}

\newcommand{\eps}{\varepsilon}

\DeclareMathOperator{\sign}{sign}

\begin{document}
\title[global existence of solutions to Euler-Poisson equations]{
Remarks on global existence of classical solution to multi-dimensional 
compressible Euler-Poisson equations with geometrical symmetry}
\author[S. Masaki]{Satoshi Masaki}
\address{Department of Mathematics\\ Kyoto University\\
Kyoto 606-8502, Japan}
\email{machack@math.kyoto-u.ac.jp}
\address{Current address: Graduate School of Information Sciences\\
Tohoku University\\ Sendai 980-8579, Japan}
\email{masaki@ims.is.tohoku.ac.jp}
\begin{abstract}
We give a necessary and sufficient condition for
the global existence of the classical solution to the Cauchy problem of the
compressible Euler-Poisson equations with radial symmetry.
We introduce a new quantity which describes the balance between the initial velocity
of the flow and the strength of the force governed by Poisson equation.
\end{abstract}
\maketitle

\section{introduction}
We consider the compressible Euler-Poisson equations:
\begin{align}
	\rho_t + \mathrm{div} (\rho v) &= 0, \label{eq:EP-1}\\
	v_t + v \cdot \nabla v &= -\l \nabla \Phi, \label{eq:EP-2}\\
	\Delta \Phi &= \rho, \label{eq:EP-3}
\end{align}
where $(t,x) \in \R_+ \times \R^n$.
These are the conservation of mass, Newton's second law, and the
Poisson equation defining, say, the electric field in terms of the charge, respectively.
The unknowns are the mass density $\rho=\rho(t,x)$,
the velocity field $v=v(t,x)$, and the potential $\Phi=\Phi(t,x)$.
$\l$ is a given physical constant.

In this paper, we assume that the unknowns have radial symmetry and
concentrate on the multi-dimensional isotropic model:
\begin{align}
	r^{n-1} \rho_t + \partial_r (r^{n-1}\rho v) = 0, \label{eq:rEP-1}\\
	v_t + v \partial_r v + \l \partial_r \Phi = 0, \label{eq:rEP-2}\\
	\partial_r (r^{n-1} \partial_r\Phi ) = r^{n-1} \rho \label{eq:rEP-3}
\end{align}
for $(t,r) \in \R_+ \times \R_+$ with initial data
\begin{equation}\label{eq:rEP-3.5}
	(\rho,v)(0,r) = (\rho_0,v_0)(r), \quad \rho_0 \ge 0.
\end{equation}
Here, $r\ge0$ denotes the distance from the origin.
Now, the unknowns are $\rho=\rho(t,r)$, $v=v(t,r)$, and $\Phi=\Phi(t,r)$.

The Euler-Poisson equations arise in many physical problems such as fluid mechanics,
plasma physics, gaseous stars, quantum gravity and semiconductors, etc.
There is a large amount of literature available on the global behavior of Euler-Poisson
and related problem, from local existence in the small $H^s$-neighborhood of
a steady state \cite{MakPW,MU-JMKU,GaCPDE} to global existence of weak solution
with geometrical symmetry \cite{CW-CMP}.
For the two-carrier types in one dimension, see \cite{WaZAMP1}.
The relaxation limit for the weak entropy solution, consult \cite{MN-ARMA}
for isentropic case, and \cite{JR-QAM} for isothermal case.
The global existence for some large class of initial data near a steady state is
obtained by Guo \cite{GuCMP} assuming the flow is irrotational.

For isotropic model, the finite time blowup for three dimensional case with the
attractive force,  pressure, and compactly supported mass density is obtained in
\cite{MP-JJAM}, and the blowup for the repulsive case in the similar settings is
deduced in \cite{PeJJAM}. In \cite{ELT-IUMJ}, the global existence/finite-time
breakdown of the strong solution is studied from the view point of critical threshold.
They give a complete criterion in one-dimensional case without spatial symmetry
and with spatial symmetry in one and four dimension.
A sufficient condition for finite-time breakdown without spatial symmetry is
obtained in \cite{CT-CMS,Chae0803}, and the complete description of the
critical threshold phenomenon for the two-dimensional restricted Euler-Poisson
equations is given in \cite{LT-SIAM}.

In this paper, applying the method in \cite{ELT-IUMJ}, we discuss the
necessary and sufficient conditions for the global existence of the 
solution to the Euler-Poisson equations with spatial symmetry \eqref{eq:rEP-1}--\eqref{eq:rEP-3.5} in multi-dimensional case.

\subsection{The Euler-Poisson equations and the Schr\"odinger-Poisson system}
The Euler-Poisson equations \eqref{eq:EP-1}--\eqref{eq:EP-3} are related to
the Schr\"odinger-Poisson system via semiclassical limit.
The global existence of the solution is applicable to the study of the
Schr\"odinger-Poisson system. Consider the Cauchy problem
\begin{equation}\tag{SP}
\left\{
	\begin{aligned}
		i h \partial_t u^h + \frac{h^2}{2}\Delta u^h &{}= \l V_p u^h, \\
		\Delta V_p &{}= |u^h|^2,\\
		u^h(0,x) &{} = a_0(x) e^{\frac{i\phi_0}{h}},
	\end{aligned}
\right.
\end{equation}
where $h$ is a positive parameter corresponding to the scaled Planck's constant.
For simplicity, we suppose that $a_0$ and $\phi_0$ in the initial data are
independent of the parameter.
In the limiting process $h\to0$, following WKB type approximation is often considered;
\begin{equation}\label{eq:WKB}
	u^h(t,x) \sim e^{i\frac{\phi(t,x)}{h}}(a(t,x)+h a_1(t,x)+h^2 a_2(t,x) + \cdots).
\end{equation}
One way to justify this approximation is to employ Madelung's transform 
\[
	u^h(t,x) = \sqrt{\rho^h(t,x)}e^{i\frac{S^h(t,x)}{h}}
\]
and consider the quantum Euler-Poisson equations
\begin{equation}\label{eq:QEP}\tag{QEP}
\left\{
	\begin{aligned}
		\partial_t \rho^h + \mathrm{div}(\rho^h \nabla S^h)&{}=0, \\
		\partial_t \nabla S^h +(\nabla S^h \cdot \nabla )\nabla S^h + \l \nabla V_p &{}= \frac{h^2}2 \nabla \( \frac{\Delta \sqrt{\rho^h}}{\sqrt{\rho^h}}\) ,\\
		\Delta V_p &{}= \rho^h, \\
		(\rho^h(0,x),\nabla S^h(0,x)) &{} = (|a_0|^2, \nabla(\phi_0+h\arg a_0 )).
	\end{aligned}
\right.
\end{equation}
The term $(h^2/2) \nabla (\Delta\sqrt{\rho^h}/\sqrt{\rho^h})$ is called
quantum pressure. Taking $h\to0$, we obtain, at least formally,
the Euler-Poisson equations
\begin{equation}\label{eq:EP}\tag{EP}
\left\{
	\begin{aligned}
		\partial_t \rho + \mathrm{div}(\rho v)&{}=0, \\
		\partial_t v +(v \cdot \nabla )v + \l\nabla V_p &{}= 0 ,\\
		\Delta V_p &{}= \rho, \\
		(\rho(0,x),v(0,x)) &{} = (|a_0|^2, \nabla\phi_0),
	\end{aligned}
\right.
\end{equation}
where $\rho=\lim_{h\to0} \rho^h$, $v=\lim_{h\to0} \nabla S^h$.
This limit is treated in \cite{ZhSIAM} and the convergence of the quadratic quantities 
\[
	|u^h|^2 \rightharpoonup \rho, \quad h\Im (\overline{u^h}\nabla u^h)
	\rightharpoonup \rho v
\]
is proved in the sense of Radon measure with some restrictive assumptions.
Though this convergence suggests that the solution $u^h$ may have the asymptotics
$u^h=e^{iS/h}(\sqrt{\rho}+o(1))$, it is not satisfactory.

Another way to justify \eqref{eq:WKB} is to employ a modified Madelung transform
\begin{equation}\label{eq:mMT}
	u^h = a^he^{i\frac{\Psi^h}{h}}
\end{equation}
and consider the system
\begin{equation}\label{eq:QEP2}\tag{$\text{QEP}^\prime$}
\left\{
	\begin{aligned}
		\partial_t a^h + (\nabla \Psi^h \cdot \nabla) a^h + \frac12 a^h \Delta \Psi^h &{}=i\frac{h}{2} \Delta a^h, \\
		\partial_t \nabla \Psi^h +(\nabla \Psi^h \cdot \nabla )\nabla \Psi^h + \l\nabla V_p &{}= 0 ,\\
		\Delta V_p &{}= |a^h|^2, \\
		(a^h(0,x),\nabla \Psi^h(0,x)) &{} = (a_0, \nabla\phi_0).
	\end{aligned}
\right.
\end{equation}
It is essential that $a^h$ takes complex value. and so, $S^h\neq \Psi^h$, in general.
The point is that the system \eqref{eq:QEP2} can be regarded as a symmetric
hyperbolic system with semilinear perturbation. It is proven in \cite{AC-SP} that this
system is locally well-posed for $0\le h \ll1$ and the solution $(a^h,\Psi^h)$
can be expanded as
\[
	a^h = a + ha_1 + h^2 a_2 + \cdots, \quad
	\Psi^h = \Psi + h \Psi_1 + h^2 \Psi_2 + \cdots.
\]
Plugging this to \eqref{eq:mMT}, we obtain WKB type estimate
\begin{equation}\label{eq:WKB2}
	u^h(t,x) = e^{i\frac{\Psi}{h}}(\beta_0(t,x) + h\beta_1(t,x) + h\beta_2(t,x) + \cdots)
\end{equation}
with $\beta_0=a e^{i\Psi_1}$ in a time interval which is small (in general) but
independent of the parameter. This method is first applied to analytic data
(\cite{PGEP}) and to Sobolev data (\cite{Grenier98}) for certain class of
nonlinearities, and it is generalized in 
\cite{AC-BKW,CR-CMP,LZ-ARMA,AC-GP,LT-MAA,AC-SP,CM-AA}.
We also note that the approximation of the form \eqref{eq:WKB2}
leads to some ill-posedness results for the ``usual'', that is,
non-scaled nonlinear Schr\"odinger equations (\cite{AC-LR,ThJDE,CM-AA}).
One verifies that the principal part of the solution $(a,\Psi)$ solves
\begin{equation}\label{eq:EP2}
\left\{
	\begin{aligned}
		\partial_t a + (\nabla \Psi \cdot \nabla) a + \frac12 a \Delta \Psi &{}=0, \\
		\partial_t \nabla \Psi +(\nabla \Psi \cdot \nabla )\nabla \Psi + \l\nabla V_p &{}= 0 ,\\
		\Delta V_p &{}= |a|^2, \\
		(a(0,x),\nabla \Psi(0,x)) &{} = (a_0, \nabla\phi_0).
	\end{aligned}
\right.
\end{equation}
Hence, we see that $\rho:=|a|^2=|\beta_0|^2$ and $v:=\nabla \Psi$ also solves \eqref{eq:EP}.

Either way, the problem of the justification of the global estimate of the form
\eqref{eq:WKB2} is closely related to the problem of global existence of the solution
to \eqref{eq:EP}. If the solution of \eqref{eq:EP} is not global and breaks down in
finite time, it immediately implies that the WKB type estimate breaks down at the same
time. The converse is not so clear. In one-dimensional case, this limit and the large time
WKB type estimate is given in \cite{LT-MAA} using the result in \cite{ELT-IUMJ}.

\subsection{Critical thresholds}
Before stating our main result, we briefly recall the part of the result in
\cite{ELT-IUMJ}. They introduce the notion of critical thresholds and
give several sufficient conditions for global existence and finite-time breakdown
in terms of the initial velocity. We restrict our attention to the positive time $t \ge 0$.
The necessary and sufficient condition for global existence is obtained in the case
$n=1$ or $4$. For a nonnegative integer $s$, we define
\[
	D^s:=
	\begin{cases}
	C([0,\I))& \text{ if } s=0, \\
	C([0,\I))\cap C^{s}((0,\I)) & \text{ if } s>0.
	\end{cases}
\]
\begin{theorem}[Critical thresholds in 1D case, \cite{ELT-IUMJ}]\label{thm:1DCT}
Suppose $n=1$, $\l<0$, $\rho_0\in D^s$, and $v_0\in D^{s+1}$ with $v_0(0)=0$
for some positive integer $s$. Then, the classical solution to
\eqref{eq:rEP-1}--\eqref{eq:rEP-3.5} is global if and only if
\begin{equation}\label{eq:1DCT}
	v_0(R) > -\sqrt{2|\l|R\int_0^R\rho_0(s)ds}\quad\text{and}\quad
	v_0^\prime(R) > -\sqrt{2|\l|\rho_0(R)}, \quad \forall R>0,
\end{equation}
where, in both inequalities, we allow the case where the both sides equal zero.
If $\rho_0$ and $v_0$ satisfy \eqref{eq:1DCT} then
the corresponding solution of \eqref{eq:rEP-1}--\eqref{eq:rEP-3.5} satisfies
\begin{align*}
	\rho &{}\in C^2([0,\I),D^s) \cap C^\I((0,\I),D^s), \\
	v &{}\in C^1([0,\I),D^{s+1}) \cap C^\I((0,\I),D^{s+1}), \\
	\Phi &{}\in C^2([0,\I),D^{s+2}) \cap C^\I((0,\I),D^{s+2}).
\end{align*}
The solution is unique in $C^2([0,\I),D^0) \times C^1([0,\I),D^1) \times
C^2([0,\I),D^2)$ and also solves \eqref{eq:EP-1}--\eqref{eq:EP-3} in the
distribution sense.
\end{theorem}
\begin{theorem}[Critical thresholds in 4D case, \cite{ELT-IUMJ}]\label{thm:4DCT}
Suppose $n=4$, $\l<0$, $\rho_0\in D^s$, and $v_0\in D^{s+1}$ with $v_0(0)=0$
for some positive integer $s$.
Let $C(r)=v_0^2(r)+|\l|r^{-2}\int_0^r \rho_0(s)s^3ds$.
The classical solution to \eqref{eq:rEP-1}--\eqref{eq:rEP-3.5} is global if and only if
both of the following conditions hold for all $R>0$:
\begin{enumerate}
\item $v_0(R) \ge 0$ if $\int_0^R \rho_0(s)s^3ds=0$;
\item $\partial_r C (R)\ge 0$ and $v_0(R)+Rv_0^\prime(R) > -\sqrt{2R\partial_r C(R)}$;
\end{enumerate}
where, in the last inequality, we allow the case where the both sides equal zero.
If $\rho_0$ and $v_0$ satisfy the above condition then
the corresponding solution of \eqref{eq:rEP-1}--\eqref{eq:rEP-3.5} satisfies
\begin{align*}
	\rho &{}\in C^2([0,\I),D^s) \cap C^\I((0,\I),D^s), \\
	v &{}\in C^1([0,\I),D^{s+1}) \cap C^\I((0,\I),D^{s+1}), \\
	\Phi &{}\in C^2([0,\I),D^{s+2}) \cap C^\I((0,\I),D^{s+2}).
\end{align*}
The solution is unique in $C^2([0,\I),D^0) \times C^1([0,\I),D^1) \times C^2([0,\I),D^2)$
and also solves \eqref{eq:EP-1}--\eqref{eq:EP-3} in the distribution sense.
\end{theorem}
\begin{remark}\label{rmk:s=0}
The above two theorems are true also in the case $s=0$. However, in that case,
$\rho$ is not differentiable. Then, introducing a new unknown
$m(t,r):=\int_0^r \rho(t,s)s^{n-1}ds$ and replacing \eqref{eq:rEP-1} with the
equation $\partial_t m +v \partial_r m=0$,  we say the solution $(\rho,v,\Phi)$ is
``classical'' in the sense that $(m,v,\Phi)$ solves this equation and
\eqref{eq:rEP-2}--\eqref{eq:rEP-3.5} in the classical sense. This is also true for
Theorems \ref{thm:attractive}, \ref{thm:3Dl-}, and \ref{thm:2Dl-}, below.
\end{remark}
\begin{remark}
The assumption $v_0(0)=0$ is natural because when we try to reconstruct
the solution $({\bf r},{\bf v},{\bf P})$ of \eqref{eq:EP-1}--\eqref{eq:EP-3},
the velocity ${\bf v}$ should be ${\bf v}(t,x)=(x/|x|)v(t,|x|)$.
\end{remark}
\begin{remark}\label{rmk:uniqueness}
If $(\rho,v,\Phi)$ is a solution to \eqref{eq:rEP-1}--\eqref{eq:rEP-3.5},
then, for any constant $c$, $(\rho,v,\Phi+c)$ also solves 
\eqref{eq:rEP-1}--\eqref{eq:rEP-3.5}. Therefore, in above theorems,
the solution is unique under certain condition on $\Phi$, such as $\Phi(0,0)=0$.
This is valid for all results below (Theorems \ref{thm:attractive}, \ref{thm:3Dl-},
\ref{thm:2Dl-} and Corollaries \ref{cor:ex1} and \ref{cor:ex2}).
\end{remark}
One-dimensional case is so simple that everything is made explicit.
On the other hand, what is special in four-dimensional case is that we can write down
some integral quantity explicitly. Essentially, their method gives necessary and
sufficient condition for all other dimensions. However, when we try to state them in
terms of the slope of the initial velocity, complex descriptions are inevitable.
This is why they give only sufficient conditions for global existence and finite-time
breakdown. They mentioned in \cite{ELT-IUMJ} that some further tedious 
calculations may enable us to obtain a complex criterion.

\subsection{Main results}
The purpose of this paper is to perform the ``further tedious calculation'' and 
determine the necessary and sufficient condition for global existence/finite-time
breakdown in the case other than $n=1,4$.
We introduce a new quantity with which we state the necessary and sufficient condition 
for global existence (we have already used in Theorem \ref{thm:4DCT}).
For $n \ge 3$, we define
\[
	C(r):=v_0(r)^2-\frac{2\l}{(n-2)r^{n-2}}\int_0^r \rho_0(s)s^{n-1}ds.
\]
This quantity represents the balance between the initial velocity and the strength of 
the force governed by the Poisson equation. This quantity clarifies the conditions for
higher dimensions. Note that 
\[
	\partial_r C(r) = 2 v_0(r) v_0^\prime(r) - \frac{2\l r\rho_0(r)}{(n-2)}
	+ \frac{2\l}{r^{n-1}}\int_0^r \rho_0(s)s^{n-1}ds,
\]
which contains the information about $v_0^\prime$. When we restrict ourselves to
the case $v_0>0$, then the use of $\partial_r C$ does not change the representation
of the conditions so much. However, in this paper, we modify the method in
\cite{ELT-IUMJ} and allow the case $v_0\le0$. The point is that some condition is
given in terms of $C$ and independent of the sign of $v_0$.
For example, if $n\ge3$ then one sufficient condition for finite-time breakdown is that
there exists $R>0$ such that $\partial_r C(R)<0$ (Theorem \ref{thm:3Dl-}, below).
Thus, the use of $C$ makes the statement slightly clearer.

Before treating the repulsive case $\l<0$, we first illustrate the result in the attractive 
case $\l>0$. Then, it turns out that the quantity $C(r)$ is very useful for the
description of the necessary and sufficient condition for the global existence.
\begin{theorem}[Critical thresholds for attractive case]\label{thm:attractive}
Suppose $\l>0$, $n\ge1$, $\rho_0\in D^s$, and $v_0\in D^{s+1}$ with
$v_0(0)=0$ for an integer $s\ge1$.
\begin{enumerate}
\item If $n =1$ or $2$ then the solution to \eqref{eq:rEP-1}--\eqref{eq:rEP-3.5}
is global if and only if $\rho_0 (r)= 0$, $v_0(r)\ge0$, and $\partial_r v_0(r) \ge 0$ 
holds for all $r \ge 0$. In particular, if $\rho_0 \not\equiv 0$ then the solution
breaks down in finite time.
\item If $n \ge3$ then the solution is global if and only if
\[
	v_0(r) \ge 0, \quad C(r) \ge 0, \quad\text{ and } \quad\partial_r C(r) \ge 0
\]
hold for all $r \ge0$.
\end{enumerate}
If $\rho_0$ and $v_0$ satisfy the condition for global existence, then
the corresponding solution of \eqref{eq:rEP-1}--\eqref{eq:rEP-3.5} satisfies
\begin{align*}
	\rho &{}\in C^2([0,\I),D^s) \cap C^\I((0,\I),D^s), \\
	v &{}\in C^1([0,\I),D^{s+1}) \cap C^\I((0,\I),D^{s+1}), \\
	\Phi &{}\in C^2([0,\I),D^{s+2}) \cap C^\I((0,\I),D^{s+2}).
\end{align*}
The solution is unique in $C^2([0,\I),D^0) \times C^1([0,\I),D^1) \times
C^2([0,\I),D^2)$ and also solves \eqref{eq:EP-1}--\eqref{eq:EP-3}
in the distribution sense.
\end{theorem}
In this case, the global existence of the solution is completely characterized as the
non-negativity and non-decreasing property of the quantity $C$. Roughly speaking,
if $C(R)$ is negative then the attractive force is so strong that the characteristic curve
starting at $r=R$ reaches $r=0$ in finite time. If $\partial_r C(R)$ is negative then it
implies that the shock is formed because outer wave propagates slower than inner
wave does. The proof appears somewhere.

We now turn to the repulsive case $\l<0$. In spite the fact that $C$ becomes always
nonnegative, the situation becomes more complicated; negative $v_0$ is allowed,
that is, the property $v_0<0$ does not necessarily lead to finite-time breakdown.
Moreover, positive $\partial_r C$ does not necessarily gives the global existence, 
either. We introduce the notion of \emph{pointwise condition for finite-time
breakdown} (\emph{PCFB}, for short) which is a necessary and sufficient condition
for finite-time blowup given only with the information of initial data at $r=R$.
Rigorous definition is given in Definition \ref{def:PCFB}. We also introduce the
quantity $A$:
\begin{align*}
	A(r) &:=\frac{2|\l| m_0(r)}{n-2}, & C(r) &:= v_0^2(r) + A(r) r^{-(n-2)}.
\end{align*}
We now state our main result.
\begin{theorem}\label{thm:3Dl-}
Suppose $\l<0$, $n\ge3$, $\rho_0\in D^s$, and $v_0\in D^{s+1}$ with 
$v_0(0)=0$ for an integer $s\ge1$. Then, the classical solution of 
\eqref{eq:rEP-1}--\eqref{eq:rEP-3.5} breaks down in finite time if and only if
there exists $R$ such that one of the following PCFB (given in Propositions
\ref{prop:3Dl-v+}, \ref{prop:3Dl-v0}, and \ref{prop:3Dl-v-}, below) is met.
On the other hand, the classical solution is global if and only if, for all $r>0$,
the PCFB does not hold. Moreover, if the condition for global existence is satisfied,
then the corresponding solution satisfies
\begin{align*}
	\rho &{}\in C^2([0,\I),D^s) \cap C^\I((0,\I),D^s), \\
	v &{}\in C^1([0,\I),D^{s+1}) \cap C^\I((0,\I),D^{s+1}), \\
	\Phi &{}\in C^2([0,\I),D^{s+2}) \cap C^\I((0,\I),D^{s+2}).
\end{align*}
Furthermore, it is unique in $C^2([0,\I),D^0) \times C^1([0,\I),D^1) \times
C^2([0,\I),D^2)$ and also solves \eqref{eq:EP-1}--\eqref{eq:EP-3} in
the distribution sense.
\end{theorem}
As stated above, this theorem holds also for $s=0$, see Remark \ref{rmk:s=0}.
\begin{proposition}[PCFB for $v_0>0$]\label{prop:3Dl-v+}
Suppose $\l<0$, $n\ge3$, and $v_0(R)>0$.
Then, the PCFB is that either one of following three conditions holds:
\begin{enumerate}
\item $\partial_r C(R)<0$;
\item $\partial_r C(R)=0$ and
\[
	\frac{1}{v_0(R)} - \frac{\partial_r A(R)}{2} \int_R^\I
	\frac{y^{-(n-2)}}{(C(R)-A(R)y^{-(n-2)})^{3/2}}dy <0;
\]
\item $0<\partial_r C(R) < \partial_r A(R)R^{-(n-2)}$ and
\[
	\frac{1}{v_0(R)} + \frac{1}{2}
	\int_R^{(\frac{\partial_r A(R)}{\partial_r C(R)})^{\frac{1}{n-2}}}
	\frac{\partial_r C(R) - \partial_r A(R) y^{-(n-2)}}{(C(R)-A(R)y^{-(n-2)})^{3/2}}dy \le 0.
\]
\end{enumerate}
\end{proposition}
\begin{proposition}[PCFB for $v_0=0$]\label{prop:3Dl-v0}
Suppose $\l<0$, $n\ge3$, and $v_0(R)=0$.
Then, the PCFB is that either one of following three conditions holds:
\begin{enumerate}
\item $\partial_r C(R)<0$;
\item $\partial_r C(R)=0$ and
\begin{enumerate}
	\item $n=3$;
	\item $n=4$ and $v_0^\prime(R) R< 0$;
	\item $n\ge 5$ and $v_0^\prime(R) R<- \frac{(n-2)\sqrt{C(R)}}{2}(1-I_n)$,
\end{enumerate}
where $I_n$ is a constant given by
\[
	I_n:=\int_1^\I \( (1-y^{-2})^{-\frac1{n-2}} -1\) dy<1;
\]
\item $\partial_r C(R) > 0$ and
\begin{enumerate}
	\item $n=3$ and
	\begin{align*}
		v_0^\prime R \le{}& -\frac34 \sqrt{C +R\partial_r C} \\
		&{} +\frac{\sqrt{C}}{2}\(1-\frac{R\partial_r C}{2C}\)\log \(\frac{\sqrt{C}+\sqrt{C+R\partial_r C}}{\sqrt{R\partial_r C}}\);
	\end{align*}
	\item $n=4$ and $v_0^\prime R \le - \sqrt{2R \partial_r C}$;
	\item $n\ge 5$ and
	\begin{align*}
	v_0^\prime R \le{}& - \frac{(n-2)^{\frac12}(R\partial_r C)^{\frac32}}{4C}
	\(1+\frac{(n-2)C}{R\partial_r C}\)^{\frac{n}{2(n-2)}}\\
	&{}-\frac{(n-2)C^{\frac12}}{2}\(1-\frac{R\partial_r C}{2C}\)
	 \\
	&{}\times \left[ \(1+\frac{R\partial_r C}{(n-2)C}\)^{\frac12}
	-\int_{\(1+\frac{R\partial_r C}{(n-2)C}\)^{\frac12}}^\I \( (1-y^{-2})^{-\frac1{n-2}} -1\) dy\right] .
\end{align*}
\end{enumerate}
Here, we omit $R$ variable in $C$, $\partial_r C$, and $v_0^\prime$, for simplicity.
\end{enumerate}
\end{proposition}
\begin{proposition}[PCFB for $v_0<0$]\label{prop:3Dl-v-}
Suppose $\l<0$, $n\ge3$,  and $v_0(R)<0$.
Then, the PCFB is that $A(R)=0$ or either one of following five conditions holds:
\begin{enumerate}
\item $\partial_r C(R)<0$;
\item $\partial_r C(R)=0$ and
\[
	\frac{1}{|v_0(R)|} - \frac{1}{2} \int_R^\I
	\frac{\partial_r A(R) y^{-(n-2)}}{(C(R)-A(R)y^{-(n-2)})^{3/2}}dy <2\partial_r t_*(R) ;
\]
\item $0<\partial_r C(R) \le \partial_r A(R)R^{-(n-2)}$ and
\[
	\frac{1}{|v_0(R)|} + \frac{1}{2}
	\int_R^{(\frac{\partial_r A(R)}{\partial_r C(R)})^{\frac{1}{n-2}}}
	\frac{\partial_r C(R) - \partial_r A(R) y^{-(n-2)}}{(C(R)-A(R)y^{-(n-2)})^{3/2}}dy \le 2\partial_r t_*(R);
\]
\item $\partial_r A(R)R^{-(n-2)}<\partial_r C(R) < \partial_r A(R)(R^{-(n-2)}+v_0(R)^2/A(R))$ and
\[
	\frac{1}{|v_0(R)|} + \frac{1}{2}
	\int_R^{(\frac{\partial_r A(R)}{\partial_r C(R)})^{\frac{1}{n-2}}}
	\frac{\partial_r C(R) - \partial_r A(R) y^{-(n-2)}}{(C(R)-A(R)y^{-(n-2)})^{3/2}}dy \le \max(0,2\partial_r t_*(R));
\]
\item $\partial_r A(R)(R^{-(n-2)}+v_0(R)^2/A(R))\le\partial_r C(R)$,
\end{enumerate}
where
\[
	t_*(R):= \(A(R)C(R)^{-\frac{n}{2}}\)^{\frac1{n-2}} \int_1^{R\(\frac{A(R)}{C(R)}\)^{-\frac1{n-2}}}
	\frac{dz}{\sqrt{1 -  z^{-(n-2)}}}.
\]
\end{proposition}
These conditions are very complex but explicit. Once we know the initial density
$\rho_0$ and the initial velocity $v_0$, we can calculate the condition.
Of course, in the four-dimensional case, the condition obtained by Propositions
\ref{prop:3Dl-v+}, \ref{prop:3Dl-v0}, and \ref{prop:3Dl-v-} becomes
the same one as in Theorem \ref{thm:4DCT}.
\begin{corollary}\label{cor:alt4D}
If $n=4$, the PCFB given in Propositions \ref{prop:3Dl-v+}, \ref{prop:3Dl-v0},
and \ref{prop:3Dl-v-} is reduced to the following condition:
\begin{enumerate}
\item $A(R)=0$ and $v_0(R)<0$.
\item $\partial_r C(R) <0$; 
\item $\partial_r C(R) =0$ and $v_0(R) + v_0^\prime(R) R <0$;
\item $\partial_r C(R) >0$ and $v_0(R) + v_0^\prime(R) R \le - \sqrt{2R\partial_r C(R)}$.
\end{enumerate}
\end{corollary}
\smallbreak

In the two-dimensional case, the quantities $A$ and $C$ have different definitions.
We introduce
\[
	A(r) := 2 |\l| m_0(r)
\]
and
\[
	C(r) := v_0(r)^2 - A(r)\log r.
\]
With these quantities, we obtain a similar theorem.
\begin{theorem}\label{thm:2Dl-}
Suppose $\l<0$, $n=2$, and $\rho_0\in D^s$, and $v_0\in D^{s+1}$ with
$v_0(0)=0$ for an integer $s\ge1$. Then, the classical solution of 
\eqref{eq:rEP-1}--\eqref{eq:rEP-3.5} breaks down in finite time if and only if
there exists $R$ such that one of the following PCFB (given in Propositions
\ref{prop:2Dl-v+}, \ref{prop:2Dl-v0}, and \ref{prop:2Dl-v-}, below) is met.
On the other hand, the classical solution is global if and only if,  for all $r>0$,
the PCFB does not hold. If the condition for global existence is satisfied, then
the corresponding solution satisfies
\begin{align*}
	\rho &{}\in C^2([0,\I),D^s) \cap C^\I((0,\I),D^s), \\
	v &{}\in C^1([0,\I),D^{s+1}) \cap C^\I((0,\I),D^{s+1}), \\
	\Phi &{}\in C^2([0,\I),D^{s+2}) \cap C^\I((0,\I),D^{s+2}).
\end{align*}
Furthermore, it is unique in $C^2([0,\I),D^0) \times C^1([0,\I),D^1) \times
C^2([0,\I),D^2)$ and also solves \eqref{eq:EP-1}--\eqref{eq:EP-3} in
the distribution sense.
\end{theorem}
\begin{proposition}[PCFB for $v_0>0$]\label{prop:2Dl-v+}
Suppose $\l<0$, $n=2$, and $v_0(R)>0$. Then, the PCFB is that
\[
	v_0^\prime (R) < \frac{A(R)}{2 R v_0(R)} 
\]
($\Leftrightarrow \exp(-\partial_r C(R)/\partial_r A(R)) > R$)
and that either one of following conditions holds:
\begin{enumerate}
\item $\rho_0(R)=0$ ($\partial_r A(R) =0$);
\item $\partial_r A(R) >0$ and
\[
	\frac{1}{v_0(R)} + \frac{1}{2} \int_R^{\exp \( -\frac{\partial_r C(R)}{\partial_r A(R)}\)}
	\frac{\partial_r C(R)+\partial_r A(R)\log y}
	{(C(R)+A(R)\log y)^{3/2}}dy \le 0.
\]
\end{enumerate}
\end{proposition}
\begin{proposition}[PCFB for $v_0=0$]\label{prop:2Dl-v0}
Suppose $\l<0$, $n=2$, and $v_0(R)=0$. Then, the PCFB is $A(R)>0$
and that either one of following conditions holds:
\begin{enumerate}
\item $\rho_0(R)=0$ ($\partial_r A(R) =0$);
\item $\partial_r A(R) >0$ and
\begin{align*}
	R v_0^\prime(R) \le& - \frac{\sqrt{A(R)R\partial_r A(R)}}{2}
	e^{\frac{A(R)}{R\partial_r A(R)}}  \\
	&+ \frac{2A(R) - R\partial_r A(R)}{4} \int_1^{e^{\frac{A(R)}{R\partial_r A(R)}}}\frac{dz}{\sqrt{\log z}}.
\end{align*}
\end{enumerate}
\end{proposition}
\begin{proposition}[PCFB for $v_0<0$]\label{prop:2Dl-v-}
Suppose $\l<0$, $n=2$, and $v_0(R)<0$. Then, the PCFB is $A(R)=0$ or either one
of following conditions holds (we omit all $R$ variables, for simplicity):
\begin{enumerate}
\item $\rho_0=0$ ($\partial_r A=0$);
\item $\partial_r A>0$ and
\begin{enumerate}
\item $\partial_r (v_0^2) \ge A/R + (\partial_r A) \log (R e^{A/C})$;
\item $A/R \le \partial_r (v_0^2) < A/R + (\partial_r A) \log (R e^{A/C})$ and
\[
	\frac{1}{|v_0|} + \frac{1}{2} \int_R^{\exp \( -\frac{\partial_r C}{\partial_r A}\)}
	\frac{\partial_r C+\partial_r A\log y}
	{(C+A\log y)^{3/2}}dy \le \max(0,2 \partial_r t_*);
\]
\item $\partial_r (v_0^2) < A/R$ and
\[
	\frac{1}{|v_0|} + \frac{1}{2} \int_R^{\exp \( -\frac{\partial_r C}{\partial_r A}\)}
	\frac{\partial_r C+\partial_r A\log y}
	{(C+A\log y)^{3/2}}dy \le 2 \partial_r t_*,
\]
\end{enumerate}
\end{enumerate}
where
\[
	t_*=t_*(R):= \frac{R}{A(R)^{1/2}e^{v_0(R)^2/A(R)}}\int_1^{e^{v_0(R)^2/A(R)}}
	\frac{dz}{\sqrt{\log z}}.
\]
\end{proposition}

\subsection{Some applications}
\begin{example}\label{ex:ex1}
In the following cases, \eqref{eq:rEP-1}--\eqref{eq:rEP-3.5} has a unique global
solution, and the solution solves \eqref{eq:EP-1}--\eqref{eq:EP-3} in 
the distribution sense.
\begin{enumerate}
\item $n=1$, $\l<0$, and
\[
	\rho_0(r) = e^{-r}, \quad v_0(r)= \sqrt{\frac{-\l}{e^r+e^{1/r}}}\sin r.
\]
\item $n=2$, $\l<0$, and
\[
	\rho_0(r) = \frac1{1+r^2}, \quad v_0(r)= \sqrt{- \l r}.
\]
\end{enumerate}
\end{example}
\begin{corollary}\label{cor:ex1}
Let $\l>0$ or $n \ge 3$. Suppose $\rho_0\in D^0 \cap L^1((0,\I),r^{n-1}dr)$ is
not identically zero and $v_0\in D^1$ satisfies $v_0(0)=0$ and $v_0 \to 0$ as
$r\to \I$. Then, the solution of \eqref{eq:rEP-1}--\eqref{eq:rEP-3.5} is global
if and only if $\l>0$ and $n\ge 3$, and the initial data is of particular form
\[
	v_0(r) = \sqrt{\frac{2\l}{(n-2) r^{n-2}}\int_0^r \rho_0(s)s^{n-1}ds}.
\]
If $\l>0$, $n \ge 3$, $\rho_0\in D^s \cap L^1((0,\I),r^{n-1}dr)$ for $s\ge1$,
and $v_0$ is as above, then $v_0 \in D^{s+1}$ and the corresponding
solution satisfies
\begin{align*}
	\rho &{}\in C^2([0,\I),D^s) \cap C^\I((0,\I),D^s), \\
	v &{}\in C^1([0,\I),D^{s+1}) \cap C^\I((0,\I),D^{s+1}), \\
	\Phi &{}\in C^2([0,\I),D^{s+2}) \cap C^\I((0,\I),D^{s+2}).
\end{align*}
Furthermore, it is unique in $C^2([0,\I),D^0) \times C^1([0,\I),D^1) \times
 C^2([0,\I),D^2)$ and also solves \eqref{eq:EP-1}--\eqref{eq:EP-3} in
the distribution sense.
\end{corollary}
\begin{proof}
In the case where $n=1$, $2$ and $\l>0$,  we deduce from Theorem
\ref{thm:attractive} that the solution breaks down in finite time because $\rho_0$
is nontrivial. Let $n\ge3$, then the assumptions $\rho_0\in L^1((0,\I),r^{n-1}dr)$
and $v_0 \to 0$ as $r\to \I$ imply $C(R)\to 0$ as $R \to \I$.
Since $C(0)=0$ by assumption, we see from Theorems \ref{thm:attractive} and
\ref{thm:3Dl-} that the solution is global only if $C \equiv 0$.
In the case $\l<0$, $C \equiv 0$ implies $\rho_0\equiv 0$, which is excluded
by assumption. In the case $\l>0$, the solution is global if we take the positive root
$v_0(R) = \sqrt{\frac{2\l}{(n-2) R^{n-2}}\int_0^R \rho_0(s)s^{n-1}ds}$.
\end{proof}
\begin{remark}
In Corollary \ref{cor:ex1}, the global solution has an explicit representation.
This is because we can solve \eqref{eq:ODE2-1}, below, explicitly by separation of
variables. In this corollary, the case $\l<0$, $n=1$ and the case $\l<0$, $n=2$ are
excluded. If $\l<0$ and $n=2$ then it is not clear whether or not the assumption of 
Corollary \ref{cor:ex1} leads to nonexistence of global solution, but following
another non-existence result holds. On the other hand, the case where $\l<0$ and
$n=1$ must be excluded since the first example in Example \ref{ex:ex1} is a counter
example. This example also suggests that the following different version also fails
if $n=1$.
\end{remark}
\begin{corollary}
Let $\l<0$ and $n\ge2$. Suppose $\rho_0\in D^0$ is not identically zero and
$v_0\in D^1$ satisfies $v_0(0)=0$. Suppose, in addition, that there exists a
sequence $\{r_j\}_{j\ge1}$ with $r_j\to\I$ as $j\to\I$ such that $v_0(r_j)=0$ 
for all $j\ge1$, $\limsup_{j\to\I}r_jv_0^\prime(r_j)<\I$, and
$r^n_j\rho_0(r_j)\to0$ as $j\to\I$. Then, the solution of
\eqref{eq:rEP-1}--\eqref{eq:rEP-3.5} breaks down in finite time.
\end{corollary}
\begin{proof}
In the $n\ge 3$ case, $v_0(r_j)=0$ leads to
\[
	\partial_r C(r_j) = \frac{2|\l|(r_j^n\rho_0(r_j)-(n-2)\int_0^{r_j}\rho_0(s)s^{n-1}ds)}{(n-2)r_j^{n-1}}.
\]
Since $\rho_0$ is nontrivial, $\int_0^{r_j}\rho_0(s)s^{n-1}ds>0$ for large $j$. 
Moreover, $r_j^n\rho_0(r_j)\to0$ as $j\to\I$ by assumption.
Hence, we conclude that $\partial_r C(r_j)<0$ for large $j$,
which is a sufficient condition for finite-time breakdown.

Let us proceed to the two dimensional case. We now show that, if $j$ is sufficiently
large, then the PCFB for $R=r_j$ (given in Proposition \ref{prop:2Dl-v0}) is satisfied 
and so the solution breaks down in finite time. Since $\rho_0$ is nontrivial, we can 
suppose $A(r_j)=2|\l|\int_0^{r_j}\rho_0(s)sds>0$. The case $\rho_0(r_j)=0$ is
trivial and so we now suppose $\partial_r A(r_j)>0$.
It suffices to prove that the inequality
\begin{align}\label{eq:ex151}
	r_j v_0^\prime(r_j) \le& - \frac{\sqrt{A(r_j)r_j \partial_r A(r_j)}}{2}
	e^{\frac{A(r_j)}{r_j\partial_r A(r_j)}}  \\
	&+ \frac{2A(r_j) - r_j\partial_r A(r_j)}{4} \int_1^{e^{\frac{A(r_j)}{r_j\partial_r A(r_j)}}}\frac{dz}{\sqrt{\log z}} \nonumber
\end{align}
is true for some $j$. Since the left hand side is upper bounded for large $j$,
by assumption, it suffices to show that the right hand side is arbitrarily large for 
large $j$. Notice that the right hand side of \eqref{eq:ex151} can be written as 
$(A(r_j)/2)f(A(r_j)/r_j\partial_r A_0(r_j))$, where
\[
	f(x) = - \frac{1}{\sqrt{x}}e^x + \int_1^{e^x}\frac{dz}{\sqrt{\log z}} -\frac1{2x}\int_1^{e^x}\frac{dz}{\sqrt{\log z}}.
\]
Since $f(1/2)=-\sqrt{2e}$ and $f^\prime(x)=(2x^2)^{-1}\int_1^{e^x}
(\log z)^{-1/2}dz\ge (2x^{5/2})^{-1}(e^x-1)$,
we see that $f(x)\to\I$ as $x\to\I$. By assumption,
$A(r_j)/r_j\partial_r A(r_j)=\int_0^{r_j}\rho_0(s)sds/r_j\rho_0(r_j)\to\I$
as $j\to\I$. Thus, the right hand side of \eqref{eq:ex151} goes to infinity as $j\to\I$.
\end{proof}
\begin{corollary}\label{cor:ex2}
Suppose $n\ge1$, $\rho_0 \equiv \rho_c >0$ is a constant, and $v_0 \equiv 0$.
Then, the solution of \eqref{eq:rEP-1}--\eqref{eq:rEP-3.5} is global
if and only if $\l<0$.  If $\l<0$ then the corresponding solution satisfies
\begin{align*}
	\rho &{}\in C^2([0,\I),D^\I) \cap C^\I((0,\I),D^\I), \\
	v &{}\in C^1([0,\I),D^\I) \cap C^\I((0,\I),D^\I), \\
	\Phi &{}\in C^2([0,\I),D^\I) \cap C^\I((0,\I),D^\I).
\end{align*}
Furthermore, it is unique in $C^2([0,\I),D^0) \times C^1([0,\I),D^1) \times
C^2([0,\I),D^2)$ and also solves \eqref{eq:EP-1}--\eqref{eq:EP-3} in the
distribution sense.
\end{corollary}
\begin{proof}
We first consider positive $\l$ case. Since $\rho_0$ is not zero, solution breaks down
if $n=1$, $2$. In the case $n\ge3$, we have $C(R)<0$ for all $R>0$,
which immediately leads to finite time breakdown.

Let us show that the solution is global if $\l<0$. The one-dimensional case is obvious
from Theorem \ref{thm:1DCT}. In the two-dimensional case, we apply Theorem
\ref{thm:2Dl-}. The PCFB is given by Proposition \ref{prop:2Dl-v0} for all $R>0$
because $v_0 \equiv 0$.
Notice that $A(R)=|\l|\rho_c r^2>0$ and $\partial_r A(R)=2|\l|\rho_c r>0$ for all $R>0$.
Therefore, in the end, we see that the solution breaks down if and only if there exists $R_0>0$ such that
\begin{align*}
	R_0 v_0^\prime(R_0) \le& - \frac{\sqrt{A(R_0)R_0\partial_r A(R_0)}}{2}
	e^{\frac{A(R_0)}{R_0\partial_r A(R_0)}}  \\
	&+ \frac{2A(R_0) - R_0\partial_r A(R_0)}{4} \int_1^{e^{\frac{A(R_0)}{R_0\partial_r A(R_0)}}}\frac{dz}{\sqrt{\log z}}.
\end{align*}
However, the left hand side is zero, and the second term of the right hand side is also zero by the 
relation $2A(R)-R\partial_r A(R)\equiv 0$.
Since the first term in the right side is negative, such $R_0$ does not exist and so the solution to
\eqref{eq:rEP-1}--\eqref{eq:rEP-3.5} is global.

We proceed to the case $n\ge3$.
The proof is the same as in two-dimensional case.
Notice that $\partial_r C(R)=4|\l|\rho_c R/n(n-2)>0$ and so that
 Proposition \ref{prop:3Dl-v0} gives the PCFB.
In the case $n=4$, it is obvious that there does not exist $R_0$ such that
$v_0^\prime (R_0)R_0 \le -\sqrt{2R_0\partial_r C(R_0)}$.
In the cases $n=3$ and $n \ge 5$, by using the fact that $C(R)=2|\l|\rho_c R^2/n(n-2)$
and so $R\partial_r C/2C\equiv 1$,
we verify nonexistence of $R_0$ for which the PCFB holds.
\end{proof}

The rest of paper is organized as follows:
We first collect some preliminary results and illustrate the strategy for proof in Section 2.
The main issues there are a reduction of the Euler-Poisson and an introduction of the notion
of pointwise condition for finite-time breakdown.
Then, we prove our main theorems in section 3.
\section{Preliminaries}
\subsection{Reduction of the Euler-Poisson equations to an ODE for characteristic curve}
We reduce the above system \eqref{eq:rEP-1}--\eqref{eq:rEP-3.5}
by employing the characteristic curve $X$
defined by an ODE
\[
	\frac{d}{dt} X(t,R) = v(t,X(t,R)), \quad X(0,R) = R
\]
and introducing the  ``mass''
\[
	m(t,r) := \int_0^r \rho(t,s) s^{n-1} ds.
\]
Then, an integration of \eqref{eq:rEP-1} yields
\begin{equation}\label{eq:rEP-1.5}\tag{\ref{eq:rEP-1}${^\prime}$}
	\partial_t m + v \partial_r m =0,
\end{equation}
which is written as
\begin{equation}\label{eq:rEP-4}
	\frac{d}{dt} m(t,X(t,R)) = 0. 
\end{equation}
Integrating \eqref{eq:rEP-3} and combining with \eqref{eq:rEP-2},
we also have
\begin{equation}\label{eq:rEP-5}
	\frac{d^2}{dt^2} X(t,R) = \frac{d}{dt} v(t,X(t,R)) = -\frac{\l m(t,X(t,R))}{(X(t,R))^{n-1}}. 
\end{equation}
Note that \eqref{eq:rEP-4} implies that the mass is conserved along
the characteristic curve.
Thus, we get an ODE for $X$:
\begin{equation}\label{eq:ODE}
	X^{\prime\prime}(t,R) = - \frac{\l m_0(R)}{X(t,R)^{n-1}}, \quad X^\prime (0,R)=v_0(R),
	\quad X(0,R) =R,
\end{equation}
where $m_0$ is the ``initial mass''  $m_0(R)=\int_0^R\rho_0(s)s^{n-1}ds$.
This reduction is the key for our analysis.
Multiply both sides by $X^\prime$ to obtain
\begin{equation}\label{eq:ODE2-1}
	(X^{\prime}(t,R))^2 = v_0(R)^2 - \frac{2\l m_0(R)}{(n-2)R^{n-2}}+
	\frac{2\l m_0(R)}{(n-2)X(t,R)^{n-2}}
\end{equation}
if $n \ge 3$ and
\begin{equation}\label{eq:ODE2-2}
	(X^{\prime}(t,R))^2 = v_0(R)^2 + {2\l m_0(R)}\log \frac{X(t,R)}{R}
\end{equation}
if $n=2$.

We now state the result about existence of $X$.
We regard $X(t,R)$ as a function $\R_+\times \R_+ \to \R$.
For a nonnegative integer $s$, we define
\[
	D^s:=
	\begin{cases}
	C([0,\I))& \text{ if } s=0, \\
	C([0,\I))\cap C^{s}((0,\I)) & \text{ if } s>0.
	\end{cases}
\]
For nonnegative integers $s_1$, $s_2$ and intervals $I_1$, $I_2$, we define
\begin{align*}
	C^{s_1,s_2}(I_1 \times I_2) =
	\{&
	f(t,x) :I_1\times I_2 \to \R  |\partial_t^a \partial_x^b f \in C(I_1\times I_2),\\ &
	\forall a\in[0, s_1], \forall b\in[0, s_2]
	\}.
\end{align*}
\begin{proposition}[Existence and regularity of solution of ODE \eqref{eq:ODE}]\label{prop:prel1}
Suppose $n \ge1$ and $\l\in\R$.
Let $s$ be a nonnegative integer and assume
$\rho_0 \in {D}^s$ and
$v_0 \in D^{s+1}$
with $v_0(0)=0$.
Then, $m_0 \in D^{s+1}$ holds, and for any $R>0$
there exists $t(R)>0$ such that $X(t,R)$ is defined from ODE \eqref{eq:ODE}
in an interval $[0,t(R))$.
Moreover, if there exists $T>0$ such that $X(t,R)>0$ holds for all $(t,R) \in [0,T) \times (0,\I)$, then
we have
\[
	X \in  C^{2,s+1}([0,T)\times (0,\I)) \cap C^{\I,s+1}((0,T)\times (0,\I)).
\]
\end{proposition}
\subsection{Local existence of the solution of \eqref{eq:rEP-1}--\eqref{eq:rEP-3.5}}
We introduce the indicator function
\begin{equation}\label{def:Ga}
	\Gamma(t,R) := \exp\(\int_0^t \partial_r u (s,X(s,R))ds\).
\end{equation}
The interpretation of $\Gamma(t,R)$ will be clear from the following lemma.
\begin{lemma}[Lemma 5.1 in \cite{ELT-IUMJ}]\label{lem:Xtosol}
Consider the Euler-Poisson equations \eqref{eq:rEP-1}--\eqref{eq:rEP-3.5}.
Let $X$ be characteristic curve, then
\[
	\Gamma(t,R) = \partial_R X(t,R).
\]
Moreover, the solution of \eqref{eq:rEP-1}--\eqref{eq:rEP-3.5} is given by
\begin{align}
	v(t,X(t,R))&{}=\frac{d}{dt}X(t,R), \\
	\rho(t,X(t,R))&{}=\frac{R^{n-1}\rho_0(R)}{X^{n-1}\Gamma(t,R)}, \\
	\partial_r v(t,X(t,X))&{}=\frac{\partial_t \Gamma(t,R)}{\Gamma(t,R)}.
\end{align}
\end{lemma}
Using the above representation of the solution,
we deduce the following Proposition,
which plays a crucial role in our proof.
\begin{proposition}[Corollary 5.2 in \cite{ELT-IUMJ}]\label{prop:main}
The smooth solution to the radial Euler-Poisson equations \eqref{eq:rEP-1}--\eqref{eq:rEP-3.5}
 is global
if and only if $\Gamma(t,R)$ is positive for all $t\ge0$ and $R \ge 0$.
On the other hand, the smooth solution to the Euler-Poisson equations
breaks down at $t=t_c$ if and only if the following equivalent conditions
are met for some $R=R_c$:
\begin{enumerate}
\item $\int_0^{t_c} \partial_r v(\tau,X(\tau,R_c)) d\tau = -\I$;
\item $\Gamma(t_c,R_c) = 0$;
\item $\partial_R X(t_c,R_c) = 0$.
\end{enumerate}
\end{proposition}
Even if it is  possible to determine a function $X$ which solves
the ODE \eqref{eq:ODE} for large time,
we can define the solution to the Euler-Poisson equations
\eqref{eq:rEP-1}--\eqref{eq:rEP-3.5} by Lemma \ref{lem:Xtosol}
as long as $X$ and $\Gamma=\partial_R X$ are positive.
\begin{proposition}[Local existence of the solution of \eqref{eq:rEP-1}--\eqref{eq:rEP-3.5}]\label{prop:prel2}
Suppose $n \ge1$ and $\l\in\R$.
Let $s$ be a nonnegative integer and assume
$\rho_0 \in {D}^s$ and
$v_0 \in D^{s+1}$ with $v_0(0)=0$.
Let $X$ be the solution of \eqref{eq:ODE} given by Proposition \ref{prop:prel1}.
Define $\Gamma$ by \eqref{def:Ga}.
If $X(t,R)>0$ and $\Gamma(t,R)>0$ hold for all $R>0$ and $t\in[0,T)$
and if $\liminf_{R\to0}\Gamma(t,R)>0$ for $t\in[0,T)$, then $X(t,0)=0$ for $t\in[0,T)$ and
\eqref{eq:rEP-1}--\eqref{eq:rEP-3.5} has a unique solution
\begin{align*}
	\rho \in{}& C^2 ([0,T), D^{s}) \cap C^\I ((0,T), D^{s}), \\
	v \in {}&C^1 ([0,T), D^{s+1}) \cap C^\I ((0,T), D^{s+1}),\\
	\Phi \in {}&C^2([0,T), D^{s+2}) \cap C^\I((0,T), D^{s+2}).
\end{align*}
\end{proposition}
\begin{remark}
In above proposition,
if $s=0$ then $\rho$ is not spatially differentiable.
In that case, we use the mass $m$ instead of $\rho$ and consider the modified equations
\eqref{eq:rEP-1.5} and \eqref{eq:rEP-2}--\eqref{eq:rEP-3.5}.
\end{remark}
\begin{remark}
If $(\rho,v,\Phi)$ is a solution to \eqref{eq:rEP-1}--\eqref{eq:rEP-3.5},
then, for any constant $c$, $(\rho,v,\Phi+c)$ also solves \eqref{eq:rEP-1}--\eqref{eq:rEP-3.5}.
Therefore, in above theorems, the solution is unique under certain condition on $\Phi$,
such as $\Phi(0,0)=0$ (See, also Remark \ref{rmk:uniqueness}).
\end{remark}
\begin{proof}
We first show $X(t,0)=0$ for $t\in[0,T)$.
Since $\liminf_{R\to0}\Gamma(t,R)>0$,
we have $R \le C X(t,R)$ for small $R$.
Then, the fact that $m_0(R)=O(R^n)$ as $R\to 0$ gives 
\[
	|X^{\prime\prime}(t,R)| =\frac{|\l| m_0(R)}{X(t,R)^{n-1}} \le C \frac{m_0(R)}{R^{n-1}}
	=O(R)
\]
as $R\to0$.
Taking the limit $R\to0$, %in the ODE \eqref{eq:ODE},
we obtain $X^{\prime\prime}(t,0)=0$ for $t \in [0,T)$.
By $X^\prime(0,0)=v_0(0)=0$ and $X(0,0)=0$, we have $X^\prime(t,0)=X(t,0)=0$ for $t \in [0,T)$.
It gives the continuities of $X$, $X^\prime$, and $X^{\prime\prime}$
(and higher time derivatives) around $R=0$:
\begin{align*}
	X \in {}& C^{2,0}([0,T)\times [0,\I)) \cap C^{2,s+1}([0,T)\times (0,\I)) \\
	& {} \cap C^{\I,0}((0,T)\times [0,\I)) \cap C^{\I,s+1}((0,T)\times (0,\I)).
\end{align*}
Then, the existence part is an immediate consequence of Lemma \ref{lem:Xtosol}.

We prove the uniqueness.
It suffices to show in the case $s=0$.
Let $(\rho_i,v_i,\Phi_i)$ ($i=1,2$) be two solutions to 
\eqref{eq:rEP-1.5}, \eqref{eq:rEP-2}--\eqref{eq:rEP-3.5}
which satisfy 
\begin{align*}
	\rho_i&{}\in C^2([0,T),D^0 ), \\
	v_i&{}\in C^1([0,T),D^1),\\
	\Phi_i&{}\in C^2([0,T),D^2).
\end{align*}
Without loss of generality, we can suppose that $\Phi_i(0,0)=0$
since, otherwise, we only have to replace $\Phi_i(t,x)$ by $\widetilde{\Phi}_i(t,x)=\Phi_i(t,x)-\Phi_i(0,0)$.
Now, solving $\frac{d}{dt}X_i(t,R)=v_i(t,X(t,R))$, we can define the characteristic curves $X_1$ and $X_2$,
and the indicator functions $\Gamma_1$ and $\Gamma_2$.
Then, we have
\begin{align*}
	X_i &{}\in C^2([0,T),D^1),&
	\Gamma_i &{} \in C^2([0,T),C((0,\I))).
\end{align*}
Since two solutions exist until $t<T$, for all $R>0$ and $\delta>0$ there exist positive constants
$c_1=c_1(R,\delta)$ and $c_2=c_2(R,\delta)$ such that
\[
	X_i(t,R) \ge c_1>0 \quad \text{ and } \quad \Gamma_i(t,R)\ge c_2>0, \quad \forall t\in[0,T-\delta].
\]
Recall that both $X_1$ and $X_2$ solve
\[
	X^{\prime\prime}(t,R) = -\l \frac{m_0(R)}{X(t,R)^{n-1}}, \quad X^\prime(0,R)=v_0(R),
	\quad X(0,R)=R.
\]
We fix $R>0$ and $\delta>0$. 
If $n=1$ then we immediately obtain $X_1(t,R)=X_2(t,R)$ for $t\in[0,T-\delta]$.
Let us proceed to the case $n\ge2$. 
Using the fact that
\[
	\left|\frac{1}{X_1(t,R)^{n-1}} -\frac{1}{X_2(t,R)^{n-1}}\right| \le \frac{n-1}{c_1^n}|X_1(t,R)-X_2(t,R)|
\]
for all $t\in[0,T-\delta]$, and applying Gronwall's lemma to
\[
	X(t,R)=R+\int_0^t X^\prime(\tau) d\tau, \quad
	X^\prime(t,R)=v_0(R)- \int_0^t \frac{\l m_0(R)}{X(\tau,R)^{n-1}}d\tau,
\]
we deduce that $X_1^\prime(t,R)=X_2^\prime(t,R)$ and $X_1(t,R)=X_2(t,R)$ hold for $t\in [0,T-\delta]$.
Since $R>0$ is arbitrary, we also have $X_1(t,0)=X_2(t,0)$ for all $t\in [0,T-\delta]$ by continuity.
Thus, we see that $X_1(t,R)=X_2(t,R)$ for all $R \ge 0$ and $t \in [0,T)$ since $\delta>0$ is also arbitrary.
Applying Lemma \ref{lem:Xtosol}, we conclude that $\rho_1=\rho_2$,
$v_1=v_2$, and so $\Phi_1=\Phi_2$.
\end{proof}
\subsection{Pointwise condition for finite-time breakdown}
Let us proceed to the discussion on global existence.
By means of Lemma \ref{lem:Xtosol} and Proposition \ref{prop:main},
it is clear that the existence of $t_c>0$ such that $\Ga(t_c,R)=0$
implies the finite-time breakdown of the solution.
The next elementary lemma suggests that the existence of $t_c>0$ such that $X(t_c,R_c)=0$ 
with some $R_c>0$ also leads to the same situation.
\begin{lemma}\label{lem:criterion1}
Let $X$ be a characteristic curve.
If $X(t_0,R_1)=X(t_0,R_2)$ for some $t_0>0$ and $0\le R_1<R_2$,
then there exist $t \in [0,t_0]$ and $R \in [R_1,R_2]$ such that $\Gamma(t,R)=0$.
In particular, if $X(t_0,R_0)=0$ for some $t_0>0$ and $R_0>0$,
then there exist $t \le t_0$ and $R \le R_0$ such that $\Gamma(t,R)=0$.
\end{lemma}
By Proposition \ref{prop:main},
to ensure the existence of the global regular solution, it suffices to 
start with the initial data for which
\[
	X(t,R)>0,\quad \forall R>0\quad \text{and}\quad \Gamma(t,R)>0, \quad \forall  R\ge0
\]
hold for all $t>0$.
Now, we introduce the notion of \emph{pointwise condition for finite-time breakdown},
which is ``almost'' the necessary and sufficient condition for the existence of $t_c \in (0,\I)$ such that $\Ga(t_c,R) =0$:
\begin{definition}\label{def:PCFB}
For fixed $R>0$, we call a 
necessary and sufficient condition for the existence of $t_c\in(0,\I)$
such that $X(t_c,R)=0$ or $\Ga(t_c,R) =0$ hold
as a pointwise condition for finite-time breakdown.
In the case of $R=0$, we regard a necessary and sufficient condition for
the existence of $t_c \in (0,\I)$ such that $\Ga(t_c,0) =0$ as a pointwise condition for finite-time breakdown.
We denote PCFB, for short.
\end{definition}
With this notion, Propositions \ref{prop:main} and \ref{prop:prel2}
are reduced as follows:
\begin{proposition}\label{prop:main2}
The local solution to the radial Euler-Poisson equations \eqref{eq:rEP-1}--\eqref{eq:rEP-3.5} given in Proposition \ref{prop:prel2}
breaks down in finite time if and only if there exist some $R\ge0$ such that
the PCFB is met.
\end{proposition}
The meaning of ``pointwise'' will be clear with this proposition.
Therefore, specification of the PCFBs is the key for the proof of 
Theorems \ref{thm:3Dl-} and \ref{thm:2Dl-}.
\subsection{Construction of the solution of \eqref{eq:EP-1}--\eqref{eq:EP-3}}
At the end of this section, we confirm that the solution of 
\eqref{eq:rEP-1}--\eqref{eq:rEP-3.5} solves the original equation.
\begin{proposition}\label{prop:recovery}
Suppose $n \ge1$ and $\l\in\R$.
Assume $\rho_0 \in D^0$ and
$v_0 \in D^{1}$ with $v_0(0)=0$.
Let $(\rho,v,\Phi)$ be a solution to \eqref{eq:rEP-1}--\eqref{eq:rEP-3.5}
given in Proposition \ref{prop:prel2}.
Then, ${\bf r}(t,x):=\rho(t,|x|)$, ${\bf v}(t,x)=(x/|x|)v(t,|x|)$,
and ${\bf P}(t,x):=\Phi(t,|x|)$ solve the Euler-Poisson equations
\begin{align}
	{\bf r}_t + \mathrm{div} ({\bf r} {\bf v}) &= 0, \label{eq:EP-1p}\tag{\ref{eq:EP-1}${}^\prime$}\\
	{\bf v}_t + {\bf v} \cdot \nabla {\bf v} &= -\l \nabla {\bf P}, \label{eq:EP-2p}\tag{\ref{eq:EP-2}${}^\prime$}\\
	\Delta {\bf P} &= {\bf r} \label{eq:EP-3p}\tag{\ref{eq:EP-3}${}^\prime$}
\end{align}
in the distribution sense.
\end{proposition}
\begin{proof}
Suppose that the solution of \eqref{eq:rEP-1}--\eqref{eq:rEP-3.5} exists 
for $t < T$.
We take a test function $\varphi(t,x) \in C^\infty_0([0,T)\times \R^n)$.
Since $(m,v,\Phi)$ solves \eqref{eq:rEP-1.5}, \eqref{eq:rEP-2}--\eqref{eq:rEP-3.5}
in the classical sense, the triplet $({\bf r},{\bf v},{\bf P})$ 
solves the \eqref{eq:EP-1p}--\eqref{eq:EP-3p} in the distribution sense in $\{|x|\neq 0\}$.
Indeed, for all $R>\eps>0$ and fixed $0<t<T$, we have
\begin{align*}
	&\int_{\eps<|x|<R} ({\bf r}_t + \mathrm{div} ({\bf r} {\bf v}))\varphi dx
	=\int_{S^{n-1}}d\omega \int_{\eps}^R \partial_r (m_t+r\partial_r m)(r) \varphi(r,\omega)dr \\
	=&{}\int_{S^{n-1}}\left[ (m_t+r\partial_r m)(R) \varphi(R,\omega)
		-(m_t+r\partial_r m)(\eps) \varphi(\eps,\omega) \right] d\omega \\
		&{}- \int_{S^{n-1}}d\omega \int_\eps^R (m_t+r\partial_r m)(r) \partial_r \varphi(r,\omega) dr=0.
\end{align*}
Hence, we only consider the case where $([0,T)\times \{0\}) \cap\mathrm{supp} \varphi \neq \emptyset$.
Put a positive small number $\eps>0$ and set $Q(\varepsilon):=\{ x\in \R^n| \sup_{1\le i\le n} |x_i| < \varepsilon \}$ and
$\chi_\eps(x):={\bf 1}_{Q(\varepsilon)}(x)$.
Then, denoting $\iint fg \,dxdt$ by $\left<f,g\right>$, we have
\begin{equation}\label{eq:recovery1}
	\left< {\bf r}_t + \mathrm{div} ({\bf r} {\bf v}), \varphi\right>=
	\left< \chi_\eps {\bf r}_t, \varphi\right>
	+\left< \chi_\eps  \mathrm{div} ({\bf r} {\bf v}), \varphi\right>.
\end{equation}
We show that the left hand side is equal to zero.
The first term of the right hand side is bounded by 
\[	
	C \int \(\int_{Q(\eps)} {\bf r}(t,x) dx\) dt + C\int_{Q(\eps)}{\bf r}(0,x) dx.
\]
We write the inverse map of $R\mapsto X(t,R)$ by $X\mapsto R(t,X)$.
This is well-defined as a map from $\R_+$ to itself and 
the two limits $R\to 0$ and $X \to 0$ are equivalent
since $\partial_R X(t,R) >0$ for all $R >0$ and $X(t,0)=0$ as long as solution exists.
By the formula of $\rho$ given in Lemma \ref{lem:Xtosol}, it holds that
\begin{align*}
	\int_{Q(\eps)} {\bf r}(t,x) dx \le{}& \int_{|x|\le \sqrt{n}\eps} {\bf r}(t,x) dx \\
	={}&\int_{0}^{\sqrt{n}\eps} \frac{R(t,x)^{n-1}\rho_0(R(t,x))}{x^{n-1}\Gamma(t,R(t,x))} x^{n-1}dx \\
	={}& \int_{0}^{R(t,\sqrt{n}\eps)} \rho_0(r) r^{n-1} dr \to 0
\end{align*}
as $\eps \to 0$ since $\rho_0(0)<\I$ by continuity.
The second term in the right hand side of \eqref{eq:recovery1}
is written as $\sum_{i=1}^n\left< \chi_\eps  \partial_i ({\bf r} {\bf v}), \varphi\right>$.
We estimate only the case $i=1$ and $n=2$:
\begin{align*}
	\left< \chi_\eps  \partial_1 ({\bf r} {\bf v}), \varphi\right>
	={}& \int_0^T \int_{-\eps}^\eps \big(({\bf r} {\bf v} \varphi)(t,\eps,x_2)-
	({\bf r} {\bf v} \varphi)(t,-\eps,x_2)\big)dx_2dt \\
	&{}-\int_0^T \int_{Q(\eps)} ({\bf r} {\bf v} \partial_1\varphi)(t,x) dxdt.
\end{align*}
Since $v(t,0)=\frac{d}{dt}X(t,0)=0$ and $v$ is continuous, we see that
$v(t,x)$ is bounded in $[0,T]\times Q(\eps)$, and so is ${\bf v}$.
Thus, both terms of the right hand side of \eqref{eq:recovery1} tend to zero because $\int_{Q(\eps)} {\bf r}(t,x) dx\to 0$ as $\eps\to0$, as checked above.
Therefore, the left hand side of \eqref{eq:recovery1} is zero since $\eps$ is arbitrary.

Similarly, we can verify that $\left< {\bf v}_t + {\bf v} \cdot \nabla {\bf v} +\l \nabla {\bf P}, \varphi\right>=0$
and $\left< \Delta{\bf P}-{\bf r}, \varphi\right>=0$.
We only note that $r^{n-1}\partial_r \Phi = \int_0^r\rho s^{n-1} ds$ by \eqref{eq:rEP-3},
and so that 
\[
	\int_{|x|\le \eps} |\nabla {\bf P}|dx \le C\int_0^{\eps}\rho s^{n-1} ds
	\to 0
\]
as $\eps \to 0$.
\end{proof}

\section{Proof of Theorems}
\subsection{Proof of Theorems \ref{thm:1DCT} and \ref{thm:4DCT}}
Let us first introduce the proof of critical thresholds in $n=1,4$
by Engelberg, Liu, and Tadmor in \cite{ELT-IUMJ}.
What is special in these cases is that the equation \eqref{eq:ODE} can be solved explicitly. 
We begin with the one-dimensional case.
\begin{proof}[Proof of Theorem \ref{thm:1DCT}]
Integrating \eqref{eq:ODE} twice, we immediately obtain
\[
	X(t,R) = R + v_0(R) t + \frac{|\l|m_0(R)}{2}t^2
\]
and so
\[
	\Gamma(t,R) = 1 + v_0^\prime(R) t+ \frac{|\l|\rho_0(R)}{2} t^2.
\]
The solution is global if and only if these two values stay positive
for all positive time.
$X(t,R)>0$ holds for all $t>0$ if and only if $v_0(R) \ge 0$ or $v_0(R)^2-|\l|Rm_0(R)/2 <0$,
and $X(t,R)>0$ holds for all $t>0$ if and only if $v_0^\prime(R) \ge 0$ or $(v_0^\prime(R))^2-|\l|\rho_0(R)/2 <0$.
Therefore, the solution is global if and only if
\[
	v_0(R) > - \sqrt{2|\l|Rm_0(R)}, \quad \text{and}\quad v_0^\prime(R) >- \sqrt{2|\l|\rho_0(R)}
\]
for all $R>0$.
Moreover, it is easy to check that the case $v_0(R)=m_0(R)=0$ and the case $v_0^\prime(R)=\rho(R)=0$
is also admissible.
\end{proof}
Let us proceed to the four-dimensional case.
\begin{proof}[Proof of Theorem \ref{thm:4DCT}]
Plugging \eqref{eq:ODE} to \eqref{eq:ODE2-1}, we see that
\[
	(X^\prime(t,R))^2= v_0(R)^2 - \frac{\l m_0(R)}{R^2} + \frac{\l m_0(R)}{X(t,R)^2}
	= C(R) - X(t,R) X^{\prime\prime}(t,R),
\]
which implies $(X(t,R)^2)^{\prime\prime}=2C(R)$. Then, integrating twice gives
\[
	X(t,R)= \sqrt{R^2 + 2v_0(R) R t + C(R) t^2} 
\]
and so
\[
	\Gamma(t,R) = \frac{2R + 2(v_0(R)+v_0^\prime(R) R)t + \partial_r C(R) t^2}{2X(t,R)}.
\]
Since $(v_0(R) R)^2 - C(R)R^2=-|\l| m_0(R)$, $X(t,R)>0$ holds for all $t>0$ if $m_0(R)>0$ or if $m_0(R)=0$ and
$v_0(R) \ge0$.
On the other hand, $\Gamma(t,R)$ stays positive for all positive time if and only if
one of the following conditions holds
\begin{enumerate}
\item $\partial_r C(R) \ge 0$ and $v_0(R)+v_0^\prime(R) R \ge0$;
\item $\partial_r C(R) > 0$, $v_0(R)+v_0^\prime(R) R <0$, and
$(v_0(R)+v_0^\prime(R) R)^2 - 2R\partial_r C(R) <0$.
\end{enumerate}
Therefore, we obtain the stated criterion.
\end{proof}
\subsection{PCFBs for repulsive $n\ge3$ case}
We prove Theorem \ref{thm:3Dl-}.
From Propositions \ref{prop:main2} and \ref{prop:recovery}, all our task is to determine the PCFB,
that is, to show Propositions \ref{prop:3Dl-v+}, \ref{prop:3Dl-v0}, and \ref{prop:3Dl-v-}.

\begin{proof}[Proof of Proposition \ref{prop:3Dl-v+}]
We first note that, by \eqref{eq:ODE} and the assumption $\l<0$,
$X^{\prime\prime}(t,R)>0$ holds as long as $X(t,R)>0$.
Since $X^\prime(0,R)=v_0(R)>0$,
we have $X^\prime(t,R)>0$, at least for small time $t \in [0,T_0]$.
Note that $X^\prime (t,R)>0$ for $t\in[0,T_0]$ implies that, for $t\in[0,T_0]$,
$X(t,R) \ge X(0,R)=R>0$ and so $X^{\prime\prime}(t,R)>0$.
Then, it means that $X^\prime$ is also increasing for $t\in[0,T_0]$.
Thus, we can choose $T_0$ arbitrarily large, that is, $X^\prime(t,R)>0$ for all $t\ge0$.
Then, for all $t\ge0$, it follows from \eqref{eq:ODE2-1} that
\[
	\int_R^{X(t,R)} \frac{dy}{\sqrt{C(R) - A(R)y^{-(n-2)}}}=t.
\]
This identity tells us that $X(t,R) \to \I$ as $t\to\I$
(This also follows from the fact that $X^\prime(t,R) \ge X^\prime(0,R) =v_0(R)>0$).
For simplicity, we omit the $R$ variable in the following.
Differentiate with respect to $R$ to obtain
\[
	\frac{\Ga(t)}{\sqrt{C - AX(t)^{-(n-2)}}}
	-\frac{1}{v_0}
	-\frac12 \int_R^{X(t)} \frac{\partial_r C
	- \partial_r Ay^{-(n-2)}}{\(C - Ay^{-(n-2)}\)^{3/2}}dy=0.
\]
We put
\[
	B(t):=\frac{\Ga(t)}{\sqrt{C - A X^{-(n-2)}}}
	=\frac{1}{v_0}
	+\frac12 \int_R^{X(t)} \frac{\partial_r C
	- \partial_r Ay^{-(n-2)}}{\(C - Ay^{-(n-2)}\)^{3/2}}dy.
\]

Assume $\partial_r C(R)<0$. Then, since $X(t) \to \I$ as $t \to \I$, 
\[
	\frac{d}{dt}B(t)=\frac{\partial_r C
	- \partial_r AX(t)^{-(n-2)}}{2\(C - AX(t)^{-(n-2)}\)^{3/2}}X^\prime(t)
	< \frac{\partial_r C}{2C^{3/2}}v_0<0
\]
holds for sufficiently large $t$. Hence, we have $B(t)\to -\I$ as $t \to \I$, and so
there always exists a time $t_0\ge 0$ such that $B(t_0) \le 0$.
We see that $\partial_r C(R)<0$ is a sufficient condition for finite-time breakdown.

Next we assume $\partial_r C(R)=0$. Then, $B(t)$ is monotone decreasing because
\[
	\frac{d}{dt}B(t)=-\frac{
	\partial_r AX(t)^{-(n-2)}}{2\(C - AX(t)^{-(n-2)}\)^{3/2}}X^\prime(t) \le 0.
\]
Therefore, there exists a time $t_0\ge 0$ such that $B(t_0) \le 0$ if and only if
$\lim_{t\to\I}B(t)<0$ (including the case $\lim_{t\to\I}B(t)=-\I$).
This condition is equivalent to
\[
	\frac{1}{v_0}
	-\frac12 \int_R^\I \frac{ \partial_r Ay^{-(n-2)}}{\(C - Ay^{-(n-2)}\)^{3/2}}dy <0.
\]

We finally assume $\partial_r C(R)>0$.
We first consider the case $(\frac{\partial_r A}{\partial_r C})^{\frac{1}{n-2}} > R$.
Then, $B(t)$ takes it minimum at a time $t=t_1\ge0$ such that 
\[
	X(t_1,R) = \(\frac{\partial_r A}{\partial_r C}\)^{\frac{1}{n-2}}> R
\]
because $\frac{d}{dt}B(t)$ is as above and $t_1$ is the time such that
$\frac{d}{dt}B(t_1)=0$.
Therefore, there exists a time $t_0$ such that $B(t_0)\le0$ if and only if
\[
	B(t_1) = \frac{1}{v_0} + \frac{1}{2}
	\int_R^{(\frac{\partial_r A}{\partial_r C})^{\frac{1}{n-2}}}
	\frac{\partial_r C - \partial_r A y^{-(n-2)}}{(C-Ay^{-(n-2)})^{3/2}}dy \le 0.
\]
We finally consider the case $(\frac{\partial_r A}{\partial_r C})^{\frac{1}{n-2}} \le R$.
However, in this case, $B$ is monotone increasing. 
Therefore, $B \ge B(0)=1/v_0>0$ for all $t \ge 0$.
\end{proof}
\begin{remark}
The argument for above proof is essentially the same as in \cite{ELT-IUMJ}.
However, this argument is not directly applicable to the case $v_0=0$.
This is because the differentiation of $\int_R^X (C-Ay^{-(n-2)})^{-1/2}dy$
produces the term $1/v_0$.
Therefore, more delicate analysis is required if $v_0=0$.
\end{remark}

\begin{proof}[Proof of Proposition \ref{prop:3Dl-v0}]
First note that we have, at least in a small time interval, $X(t,R)>0$ because $X(0,R)=R>0$.
Since $X^{\prime\prime}(t,R) > 0$ holds as long as $X(t,R)>0$ by \eqref{eq:ODE},
we can find a time $t_0>0$ such that $X^\prime(t_0,R) > X^\prime(0,R)=v_0(R)=0$.
Note that $t_0$ can be chosen arbitrarily small.
Then, repeating the argument as in the proof of Proposition \ref{prop:3Dl-v+},
we see that, $X^\prime(t,R)\ge X^\prime(t_0,R) > 0$ for all $t \ge t_0$,
which shows $X^\prime(t,R)>0$ for all $t>0$ and $X(t,R)\to \I$ as $t\to\I$.
Moreover, $X(t,R) \sim C(R)^{1/2}t$ for sufficiently large $t$
since $X^\prime(t,R) \to C(R)^{1/2}$ as $t\to\I$.
It reveals that if $\partial_r C(R)<0$ then the characteristic curves must cross
and so the solution breaks down in finite time by Lemma \ref{lem:criterion1}.

We now suppose $\partial_r C(R) \ge 0$.
We omit $R$ variable in the following.
Since $X^\prime(t) \ge 0$ for all $t\ge0$,
an integration of \eqref{eq:ODE2-1} gives
\[
	\int_R^{X(t)} \frac{dy}{\sqrt{C-Ay^{-(n-2)}}} = t.
\]
By a change of variable $z=y/R$, the left hand side is equal to 
\[
	\int_1^{X(t)/R} \frac{R dz}{\sqrt{C-AR^{-(n-2)}z^{-(n-2)}}}.
\]
We temporally assume that $v_0>0$ and take the limit $v_0 \downarrow 0$ later.
This computation is justified, for example,
by replacing $v_0$ by $X^\prime(\eps R,R)>0$ with small $\eps>0$ and taking the limit $\eps \downarrow 0$.
Differentiation with respect $R$ yields
\begin{align*}
	0=&\frac{R \partial_R(X(t)/R)}{\sqrt{C-AX(t)^{-(n-2)}}}
	+ \int_1^{X(t)/R} \frac{dz}{\sqrt{C-AR^{-(n-2)}z^{-(n-2)}}}\\
	& - R\int_1^{X(t)/R} \frac{\partial_r C-(\partial_r AR^{-(n-2)}-(n-2)AR^{-(n-1)}
	)z^{-(n-2)}}{2\(C-AR^{-(n-2)}z^{-(n-2)}\)^{3/2}}dz.
\end{align*}
For simplicity, we omit $t$ variable in $X$ and $\partial_R X$ for a while
because the following computations do not include any differentiation.
An elementary calculation shows
\begin{align}\label{eq:predef3Dl-v0}
	0
	=&\frac{\partial_RX}{\sqrt{C-AX^{-(n-2)}}}
	- \frac{X}{R\sqrt{C-AX^{-(n-2)}}}
	+ \int_1^{X/R} \frac{dz}{\sqrt{C-AR^{-(n-2)}z^{-(n-2)}}} \\
	&- \frac{R\partial_r C}{2C}
	\int_1^{X/R} \frac{C-AR^{-(n-2)}z^{-(n-2)}}{\(C-AR^{-(n-2)}z^{-(n-2)}\)^{3/2}}dz 
	\nonumber \\
	&- \frac{R\partial_r C}{2C}
	\int_1^{X/R} \frac{AR^{-(n-2)}z^{-(n-2)}}{\(C-AR^{-(n-2)}z^{-(n-2)}\)^{3/2}}dz 
	\nonumber \\
	&+ \int_1^{X/R} \frac{R(\partial_r AR^{-(n-2)}-(n-2)AR^{-(n-1)}
	)z^{-(n-2)}}{2\(C-AR^{-(n-2)}z^{-(n-2)}\)^{3/2}}dz 
	\nonumber \\
	=&\frac{\partial_RX}{\sqrt{C-AX^{-(n-2)}}}
	- \frac{X}{R\sqrt{C-AX^{-(n-2)}}}
	+ \int_1^{X/R} \frac{dz}{\sqrt{C-AR^{-(n-2)}z^{-(n-2)}}}
	\nonumber\\
	&- \frac{R\partial_r C}{2C}
	\int_1^{X/R} \frac{dz}{\(C-AR^{-(n-2)}z^{-(n-2)}\)^{1/2}}
	\nonumber \\
	&+ \frac{(-\partial_r C AR+C\partial_r AR-(n-2)AC
	)}{2CR^{n-2}}\int_1^{X/R} \frac{z^{-(n-2)}}{\(C-AR^{-(n-2)}z^{-(n-2)}\)^{3/2}}dz. \nonumber
\end{align}
It also holds that
\begin{align*}
	\frac{-\partial_r C AR+C\partial_r AR-(n-2)AC}{2CR^{n-2}}
	=\(-\frac{v_0^\prime A}{CR^{n-3}} +
	\frac{\partial_r A R -(n-2) A}{2CR^{n-2}}v_0\)v_0.
\end{align*}
Now, let us show that
\begin{equation}\label{eq:limit3Dl-v0}
	\lim_{v_0 \downarrow 0} v_0 \int_1^{X/R} \frac{z^{-(n-2)}}{\(C-AR^{-(n-2)}z^{-(n-2)}\)^{3/2}}dz =\frac{2}{AR^{-(n-2)}(n-2)}.
\end{equation}
Fix a small $\eps>0$.
Then, we have
\[
	\lim_{v_0 \downarrow 0} v_0 \int_{1+\eps}^{X/R}
	\frac{z^{-(n-2)}}{\(C-AR^{-(n-2)}z^{-(n-2)}\)^{3/2}}dz =0,
\]
since the integral is uniformly bounded with respect to $v_0$.
Moreover,
\begin{align*}
	&v_0 \int_1^{1+\eps}
	\frac{z^{-(n-2)}}{\(C-AR^{-(n-2)}z^{-(n-2)}\)^{3/2}}dz \\
	& \le \frac{2v_0(1+\eps)}{AR^{-(n-2)}(n-2)}\int_1^{1+\eps}
	\frac{AR^{-(n-2)}(n-2)z^{-(n-1)}}{2\(C-AR^{-(n-2)}z^{-(n-2)}\)^{3/2}}dz \\
	& \le \frac{2v_0(1+\eps)}{AR^{-(n-2)}(n-2)}
	\left[\(C-AR^{-(n-2)}\)^{-\frac12}
	- \(C-AR^{-(n-2)}(1+\eps)^{-(n-2)}\)^{-\frac12}\right] \\
	& \to \frac{2(1+\eps)}{AR^{-(n-2)}(n-2)}
\end{align*}
as $v_0\to 0$.
Similarly, 
\begin{align*}
	&v_0 \int_1^{1+\eps}
	\frac{z^{-(n-2)}}{\(C-AR^{-(n-2)}z^{-(n-2)}\)^{3/2}}dz \\
	& \ge \frac{2v_0}{AR^{-(n-2)}(n-2)}\int_1^{1+\eps}
	\frac{AR^{-(n-2)}(n-2)z^{-(n-1)}}{2\(C-AR^{-(n-2)}z^{-(n-2)}\)^{3/2}}dz \\
	& \to \frac{2}{AR^{-(n-2)}(n-2)}
\end{align*}
as $v_0\to 0$.
Since $\eps>0$ is arbitrary, we obtain \eqref{eq:limit3Dl-v0}.

Taking the limit $v_0 \downarrow 0$ in \eqref{eq:predef3Dl-v0},
\begin{align*}
	0=&\frac{\partial_RX}{C^{1/2}\sqrt{1-(R/X)^{n-2}}}
	- \frac{X}{RC^{1/2}\sqrt{1-(R/X)^{n-2}}}
	+ \int_1^{X/R} \frac{dz}{C^{1/2}\sqrt{1-z^{-(n-2)}}} \\
	&- \frac{R\partial_r C}{2C^{3/2}}
	\int_1^{X/R} \frac{dz}{\sqrt{1-z^{-(n-2)}}} -\frac{2v_0^\prime R}{(n-2)C}.
\end{align*}
Thus, we have
\begin{align*}
	\frac{\partial_R X(t)}{\sqrt{1-(R/X(t))^{n-2}}} =& \frac{X(t)}{R\sqrt{1-(R/X(t))^{n-2}}}
	- \int_1^{X(t)/R} \frac{dz}{\sqrt{1-z^{-(n-2)}}}\\
	&{} +\frac{R\partial_r C }{2C}\int_1^{X(t)/R} \frac{dz}{\sqrt{1-z^{-(n-2)}}} 
	 + \frac{2v_0^\prime R}{(n-2)C^{1/2}}.
\end{align*}
We denote this by $B(t)$.
\smallbreak

{\bf Case 1.} We first assume that $\partial_r C(R) =0$.
We put 
\[
	G(s) := \frac{s}{\sqrt{1-s^{-(n-2)}}} - \int_1^{s} \frac{dz}{\sqrt{1-z^{-(n-2)}}}.
\]
An elementary calculation shows, for $s>1$,
\[
	G^\prime(s) = - \frac{(n-2)s^{-(n-2)}}{2(1-s^{-(n-2)})^{3/2}} < 0,
\]
and so $G$ is monotone decreasing.
Moreover, considering the inverse map of $z \mapsto (1-z^{-(n-2)})^{-1/2}$,
we have
\[
	\int_1^{s} \frac{dz}{\sqrt{1-z^{-(n-2)}}} = 
	\frac{s-1}{\sqrt{1-s^{-(n-2)}}} + \int_{(1-s^{-(n-2)})^{-\frac12}}^\I
	\( (1-y^{-2})^{-\frac1{n-2}} -1\) dy.
\]
Therefore,
\[
	G(s) = \frac1{\sqrt{1-s^{-(n-2)}}} - \int_{(1-s^{-(n-2)})^{-\frac12}}^\I
	\( (1-y^{-2})^{-\frac1{n-2}} -1\) dy.
\]
One verifies that if $n=3$ then $\lim_{s\to\I}G(s)=-\I$.
We now put, for $n \ge 4$,
\[
	I_n:=\int_1^\I \( (1-y^{-2})^{-\frac1{n-2}} -1\) dy.
\]
For any $m> l \ge 4$ and $y\in(1,\I)$, it holds that 
\[
	(1-y^{-2})^{-\frac1{m-2}}
	< (1-y^{-2})^{-\frac1{l-2}}.
\]
This gives $I_m<I_l$ for $m>l\ge 4$.
If $n=4$ then
\begin{align*}
	I_4 ={}& \lim_{N\to\I}\int_1^N
	\( (1-y^{-2})^{-\frac12} -1\) dy \\
	={}&  \lim_{N\to\I}
	\( (N^2-1)^{\frac12}-(N -1)\) 
	= 1.
\end{align*}
Thus, we obtain
\[
	\lim_{s\to\I}G(s)={}
	\begin{cases}
		-\I & \text{ if }n=3, \\
		0 & \text{ if }n=4, \\
		1-I_n>0 &\text{ if }n\ge5.
	\end{cases}
\]
Since
\[
	B(t) =G(X(t)/R) - \frac{2v_0^\prime R}{(n-2)C^{1/2}},
\]
we conclude that there exists $t_0\in[0,\I)$ such that $\Ga(t_0) \le 0$
if and only if
\begin{enumerate}
\item $n=3$;
\item $n=4$ and $v_0^\prime R<0$;
\item $n\ge5$ and $v_0^\prime R < - \frac{(n-2)\sqrt{C}}{2}(1-I_n)$.
\end{enumerate}
\smallbreak

{\bf Case 2.} We assume that $\partial_r C(R) >0$.
We write $B(t)=H(X(t)/R)$.
Then, it holds that
\[
	\frac{d}{ds}H(s) = -\frac{(n-2)s^{-(n-2)}}{2(1-s^{-(n-2)})^{3/2}}
	+\frac{R\partial_r C}{2C(1-s^{-(n-2)})^{1/2}}.
\]
Therefore, the minimum of $H$, hence of $B$, is
\[
	H\( \(1+\frac{(n-2)C}{R\partial_r C}\)^{\frac{1}{n-2}}\).
\]
The solution breaks down in finite time if and only if
this value is less than or equal to zero.
This leads us to the condition
\begin{align*}
	v_0^\prime R \le{}&- \frac{\sqrt{(n-2)R\partial_r C}}{2}
	\(1+\frac{(n-2)C}{R\partial_r C}\)^{\frac{n}{2(n-2)}}\\
	&{} - \frac{n-2}{2}C^{\frac12}\(\frac{R\partial_r C}{2C}-1\)
	\int_1^{\(1+\frac{(n-2)C}{R\partial_r C}\)^{\frac{1}{n-2}}} \frac{dz}{\sqrt{1-z^{-(n-2)}}}.
\end{align*}
Using the identity
\begin{align*}
	\int_1^{\(1+\frac{(n-2)C}{R\partial_r C}\)^{\frac{1}{n-2}}} \frac{dz}{\sqrt{1-z^{-(n-2)}}}
	={}& \({\(1+\frac{(n-2)C}{R\partial_r C}\)^{\frac{1}{n-2}}}-1\)
	\(1+\frac{R\partial_r C}{(n-2)C}\)^{\frac12} \\
	&{} + \int_{\(1+\frac{R\partial_r C}{(n-2)C}\)^{\frac12}}^\I \( (1-y^{-2})^{-\frac1{n-2}} -1\) dy,
\end{align*}
we obtain the equivalent condition
\begin{align*}
	v_0^\prime R \le{}& -\frac{(n-2)^{\frac12}(R\partial_r C)^{\frac32}}{4C}
	\(1+\frac{(n-2)C}{R\partial_r C}\)^{\frac{n}{2(n-2)}}\\
	&{}-\frac{(n-2)C^{\frac12}}{2}\(1-\frac{R\partial_r C}{2C}\) \\
	&{}\times \left[ \(1+\frac{R\partial_r C}{(n-2)C}\)^{\frac12}
	-\int_{\(1+\frac{R\partial_r C}{(n-2)C}\)^{\frac12}}^\I \( (1-y^{-2})^{-\frac1{n-2}} -1\) dy\right].
\end{align*}
In particular, if $n=3$ or $4$, then the above integral is computable,
and we have more explicit condition
\[
	v_0^\prime R \le - \frac34 \sqrt{C +R\partial_r C}
	+\frac{\sqrt{C}}{2}\(1-\frac{R\partial_r C}{2C}\)\log \(\frac{\sqrt{C}+\sqrt{C+R\partial_r C}}{\sqrt{R\partial_r C}}\)
\]
if $n=3$ and
\[
	v_0^\prime R \le - \sqrt{2R \partial_r C}
\]
if $n=4$.
\end{proof}
\begin{proof}[Proof of Proposition \ref{prop:3Dl-v-}]
We first note that if $A(R)=0$, then $X^\prime(t,R)\equiv v_0(R)<0$.
Therefore, the solution breaks down no latter than $t=R/|v_0(R)|$ by Lemma \ref{lem:criterion1}.
Hence, we assume $A(R)>0$.
Then, since  $X^\prime(0,R)=v_0(R)<0$, $X^\prime(t,R)=- \sqrt{C - A X(t)^{-(n-2)}}$ as long as
$X^\prime(t,R) \le 0$.
Take
\begin{align*}
	t_* =& \int_{\(\frac{A}{C}\)^{\frac1{n-2}}}^R \frac{dy}{\sqrt{C - A y^{-(n-2)}}} \\
	=& \(AC^{-\frac{n}{2}}\)^{\frac1{n-2}} \int_1^{R\(\frac{A}{C}\)^{-\frac1{n-2}}}
	\frac{dz}{\sqrt{1 -  z^{-(n-2)}}}.
\end{align*}
We see that, for all $t\in[0,t_*)$, $X(t,R)>X(t_*,R) = (A(R)/C(R))^{1/(n-2)}>0$
and  $X^\prime(t,R)<X^\prime(t_*,R)=0$.
Since $X^{\prime\prime}(t_*,R)>0$ by \eqref{eq:ODE}, using the same argument as in the proof of 
Proposition \ref{prop:3Dl-v0}, we have $X^\prime(t,R)\ge0$ for all $t\ge t_*$ and so
\[
	X^\prime(t,R) =
	\begin{cases}
	- \sqrt{C(R) - A(R) X(t,R)^{-(n-2)}}, & \text{ for }t\le t_*, \\
	\sqrt{C(R) - A(R) X(t,R)^{-(n-2)}}, & \text{ for }t\ge t_*.
	\end{cases}
\]
We also obtain $X(t,R)\to\I$ as $t\to\I$. 
In the following, we omit $R$ variable. 
For sufficient large $t$, $X(t) \sim C^{1/2}t$ holds
since $X^\prime(t) \to C^{1/2}$ as  $t\to\I$.
It implies that if $\partial_r C(R)<0$ then the characteristic curves must cross
and so the solution breaks down in finite time by Lemma \ref{lem:criterion1}.
Differentiation of $X(t_*,R) = (A/C)^{1/(n-2)}$
with respect to $R$ gives
\[
	\partial_r t_* X^\prime(t_*,R) + \partial_R X(t_*,R) =
	\partial_r \(\frac{A}{C}\)^{\frac1{n-2}}.
\]
Using the fact that $X^\prime(t_*)=0$, we obtain
\[
	\partial_R X(t_*,R) = \partial_r \(\frac{A}{C}\)^{\frac1{n-2}}.
\]
Hence, if $\partial_r (A/C)^{1/(n-2)} \le 0$ then the solution 
breaks down no latter than $t_*$.

Thus, we assume $\partial_r C(R) \ge 0$ and $(\partial_r (A/C)(R))^{1/(n-2)} > 0$
in the following.
Notice that the latter condition is equivalent to the following two conditions:
\[
	\partial_r C < \partial_r A (R^{-(n-2)} + v_0^2/A), \quad
	\(\frac{A}{C}\)^{\frac1{n-2}}
	< \(\frac{\partial_r A}{\partial_r C}\)^{\frac1{n-2}}.
\]
\smallbreak

{\bf Step 1.}
We determine the condition that solution can be extended to time $t=t_*$.
For $t \le t_*$, we have
\[
	\int_{X(t)}^R \frac{dy}{\sqrt{C-Ay^{-(n-2)}}} = t.
\]
Differentiation with respect to $R$ yields
\[
	\frac{1}{\sqrt{C-AR^{-(n-2)}}} - 
	\frac{\Ga(t)}{\sqrt{C-AX(t)^{-(n-2)}}}
	- \frac12 \int^R_{X(t)}
	\frac{\partial_r C - \partial_r A y^{-(n-2)}}{(C-Ay^{-(n-2)})^{3/2}}dy = 0.
\]
For $0 \le t < t_*$, it holds that
\[
	0 < \sqrt{C-AX(t)^{-(n-2)}} \le \sqrt{C-AR^{-(n-2)}}=|v_0|.
\]
Therefore, $\Ga(t)$ has the same sign as
\[
	B_1(t):= \frac{\Ga(t)}{\sqrt{C-AX^{-(n-2)}}}
	=\frac{1}{|v_0|} - \frac12 \int^R_{X(t)}
	\frac{\partial_r C - \partial_r A y^{-(n-2)}}{(C-Ay^{-(n-2)})^{3/2}}dy.
\]
Taking time derivative, one verifies that $B_1$ takes it minimum
at $t=t_1 \in [0,t_*)$ such that 
\[
	X(t_1,R) = \min \(R, \(\frac{\partial_r A}{\partial_r C}\)^{\frac1{n-2}}\).
\]
Note that $(A/C)^{1/(n-2)} < X(t_1)$ by assumption,
and that $({\partial_r A}/{\partial_r C})^{1/(n-2)}<R$ is equivalent to
$\partial_r C > \partial_r A R^{-(n-2)}$.
Since we have already known that $\Ga(0) =1>0$,
the solution can be extended to the time $t=t_*$ UNLESS
$\partial_r C > \partial_r A R^{-(n-2)}$ and
\[
	B_1(t_1)=
	\frac{1}{|v_0|} - \frac12 \int^R_{\(\frac{\partial_r A}{\partial_r C}\)^{\frac1{n-2}}}
	\frac{\partial_r C - \partial_r A y^{-(n-2)}}{(C-Ay^{-(n-2)})^{3/2}}dy
	\le 0
\]
is satisfied. Notice that this condition is a sufficient condition for finite-time breakdown.
\smallbreak

{\bf Step 2.}
We consider the condition that the solution
can be extended from the time $t=t_*$ to $t=\I$.
For simplicity, we suppose that solutions are extended to time $t=t_*$
(we keep assuming $0\le \partial_r C < \partial_r A (R^{-(n-2)} + v_0^2/A)$
holds).
Recall that, for $t\ge t_* $, $X^\prime(t) = \sqrt{C-AX(t)^{-(n-2)}}\ge0$.
As in the case $v_0=0$, this inequality with $X^{\prime\prime}(t)>0$ gives
$X(t)\sim C^{1/2}t \to\I$ as $t \to \I$.

We define $t_{**}$ as the time that $t_{**}>t_{*}$ and $X(t_{**})=R$.
Then, we have
\[
	t_{**}-t_* = \int_{\(\frac{A}{C}\)^{\frac1{n-2}}}^R \frac{dy}{\sqrt{C - A y^{-(n-2)}}} = t_*.
\]
Therefore, $t_{**}=2t_*$ and
\[
	\int_R^{X(t)} \frac{dy}{\sqrt{C - A y^{-(n-2)}}} = t-2t_*
\]
for all $t \ge t_*$.
As in the previous step, we set
\[
	B_2(t):= \frac{\Ga(t)}{\sqrt{C-AX(t)^{-(n-2)}}}
	=\frac{1}{|v_0|} + \frac12 \int_R^{X(t)}
	\frac{\partial_r C - \partial_r A y^{-(n-2)}}{(C-Ay^{-(n-2)})^{3/2}}dy
	-2\partial_r t_*.
\]
$B_2(t)$ and $\Ga(t)$ has the same sign for $t \ge t_*$.
We also note that $B_2(t)\to \I$ as $t \downarrow t_*$ because
$\Ga(t_*)>0$ and $\sqrt{C-AX(t)^{-(n-2)}} \to 0$ as $t \downarrow t_*$.
It holds that
\[
	\frac{d}{dt} B_2(t) =
	\frac{\partial_r C - \partial_r A X(t)^{-(n-2)}}{2(C-AX(t)^{-(n-2)})^{3/2}}X^\prime(t).
\]

\begin{enumerate}
\item If $\partial_r C(R)=0$ then $B_2$ is monotone decreasing because $\frac{d}{dt}B_2(t) \le 0$.
Therefore, solution can be extended to $t=\I$ if and only if
\[
	\lim_{t\to\I}B_2(t)=\frac{1}{|v_0|} - \frac12 \int_R^\I
	\frac{\partial_r A y^{-(n-2)}}{(C-Ay^{-(n-2)})^{3/2}}dy
	-2\partial_r t_* \ge 0.
\]
\item If $\partial_r C(R) > 0$ then $B_2$ takes it minimum at
$t=t_2$ such that $X(t_2)=({\partial_r A}/{\partial_r C})^{1/(n-2)}$.
Therefore, solution can be extended to $t=\I$ if and only if
\[
	B_2(t_2)=\frac{1}{|v_0|} + \frac12 \int_R^{\(\frac{\partial_r A}{\partial_r C}\)^{\frac1{n-2}}}
	\frac{\partial_r C - \partial_r A y^{-(n-2)}}{(C-Ay^{-(n-2)})^{3/2}}dy
	-2\partial_r t_* >0 .
\]
\end{enumerate}
\end{proof}
Before proceeding to the two-dimensional case, let us see that
Theorem \ref{thm:3Dl-} gives the same criterion as in Theorem \ref{thm:4DCT}
if $n=4$.
Namely, we prove Corollary \ref{cor:alt4D}.

\begin{proof}[Proof of Corollary \ref{cor:alt4D}]
Before the proof, we prepare some elementary computations.
We note that
\begin{align*}
	&\int_R^{\sqrt{\frac{\partial_r A}{\partial_r C}}}
	\frac{\partial_r C - \partial_r A y^{-2}}{(C-Ay^{-2})^{3/2}}dy \\
	&{} = \frac{\partial_r C}{C}\int_R^{\sqrt{\frac{\partial_r A}{\partial_r C}}}
	\frac{y}{(Cy^2-A)^{1/2}}dy
	 + \frac{A\partial_r C-C \partial_r A}{C} \int_R^{\sqrt{\frac{\partial_r A}{\partial_r C}}}
	\frac{y}{(Cy^2-A)^{3/2}}dy \\
	&{} = \frac{\partial_r C}{C^2}\left[ \(C \frac{\partial_r A}{\partial_r C}-A \)^{\frac12} - (CR^2-A)^{\frac12}\right]\\
	&\quad + \frac{A\partial_r C-C \partial_r A}{C^2}
	\left[ (CR^2-A)^{-\frac12}-\(C \frac{\partial_r A}{\partial_r C}-A \)^{-\frac12}\right]\\
	&=\frac{2(\partial_r C)^{\frac12}}{C^2}(C\partial_r A-A \partial_r C)^{\frac12}
	- \frac{|v_0|R}{C^2}\partial_r C - \frac{C\partial_r A-A \partial_r C}{C^2 |v_0|R},
\end{align*}
and that
\begin{align*}
	&\frac{1}{|v_0|} - \frac{|v_0|R}{2C^2}\partial_r C - \frac{C\partial_r A-A \partial_r C}{2C^2 |v_0|R}\\
	&{} = \(\frac{1}{|v_0|} + \frac{v_0^2R^2 +A}{2C^2 |v_0|R}\partial_r C
	 -\frac{\partial_r A}{2C |v_0|R}\)-\frac{|v_0|R}{C^2}\partial_r C \\
	&{} = \frac{2CR + R^2\partial_r C-\partial_r A}{2C |v_0|R} -\frac{|v_0|R}{C^2}\partial_r C \\
	&{} = (\sign v_0)\(\frac{v_0 + R v_0^\prime}{C} -\frac{v_0R}{C^2}\partial_r C\) \\
	&{} = (\sign v_0)\partial_r \(\frac{v_0R}{C}\),
\end{align*}
where we have used $v_0^2 R^2 + A =CR^2$ and
\begin{align*}
	&2CR + R^2\partial_r C-\partial_r A \\
	&= \(2v_0^2 R + \frac{2\l m_0}{R}\)
	+ \(2v_0 v_0^\prime R^2 + \partial_r A - \frac{2\l m_0}{R}\)
	-\partial_r A \\
	&= 2v_0 R(v_0 + R v_0^\prime).
\end{align*}
It also holds that
\begin{align*}
	t_* &{}= \sqrt{AC^{-2}} \int_1^{R\(\frac{A}{C}\)^{-\frac12}}
	\frac{dz}{\sqrt{1 -  z^{-2}}}
	= \sqrt{AC^{-2}}\int_1^{R^2\(\frac{C}{A}\)}\frac{dz}{2\sqrt{z -  1}}\\
	&{}= \sqrt{AC^{-2}} \sqrt{\frac{R^2 C}{A}-1} = \frac{|v_0| R}{C}.
\end{align*}

From Propositions \ref{prop:3Dl-v+}, \ref{prop:3Dl-v0}, and \ref{prop:3Dl-v-},
we see that $\partial_r C<0$ is the sufficient condition for blow-up.
Moreover, the PCFB in the case $\partial_r C=0$ is
\[
	\frac{(\partial_r C)^{\frac12}}{C^2}(C\partial_r A-A \partial_r C)^{\frac12}
	+ \partial_r \(\frac{v_0R}{C}\) = \frac{v_0 + Rv_0^\prime}{C}<0
\]
if $v_0 >0$, 
\[
	R v_0^\prime  <0
\]
if $v_0=0$, and
\[
	\frac{(\partial_r C)^{\frac12}}{C^2}(C\partial_r A-A \partial_r C)^{\frac12}
	- \partial_r \(\frac{v_0R}{C}\)+2\partial_r \(\frac{v_0R}{C}\) = \frac{v_0 + Rv_0^\prime}{C}<0
\]
if $v_0<0$.
Hence, the PCFB is summarized as 
$v_0 + Rv_0^\prime < 0$.

Let us proceed to the case $\partial_r C > 0$.
If $v_0>0$ then Proposition \ref{prop:3Dl-v+} implies that the PCFB is
$\partial_r C < \partial_r A R^{-2} \Leftrightarrow v_0+Rv_0^\prime < C/v_0$
and
\begin{equation}\label{eq:4DCTalt1}
	\frac{(\partial_r C)^{\frac12}}{C^2}(C\partial_r A-A \partial_r C)^{\frac12}
	+ \partial_r \(\frac{v_0R}{C}\)\le 0.
\end{equation}
We put $\alpha=v_0+Rv_0^\prime$, $\beta=v_0R\partial_r C/C>0$,
and $\gamma=\partial_r C(C\partial_r A-A \partial_r C)/C^2$.
Note that, by assumption, we have
$0<A/C<R^2< \partial_r A/\partial_r C$,
which implies $\gamma>0$.
Then, \eqref{eq:4DCTalt1} can be written as
$\alpha \le \beta - \sqrt{\gamma}$.
We make this condition clearer.
An elementary computation shows that
$\delta:=\gamma+2\alpha\beta-\beta^2=2 R \partial_r C >0$,
and that
$\delta-2\alpha\beta=\frac{R\partial_r C(-\partial_r C + \partial_r AR^{-2})}{C}>0$.
The latter one means $\beta^2< \gamma$.
Thus, the inequality $\alpha\le \beta - \sqrt{\gamma}<0$
is reduced to $\alpha\le-\sqrt{\gamma+2\alpha\beta-\beta^2}=-\sqrt{\delta}$,
that is,
$v_0 + R v_0^\prime \le -\sqrt{2R\partial_r C}$.
This condition is stronger than $\partial_r C < \partial_r A R^{-2} \Leftrightarrow v_0+Rv_0^\prime < C/v_0$.

If $\partial_r C>0$ and $v_0=0$, then it immediately follows from Proposition \ref{prop:3Dl-v0} that $\alpha \le -\sqrt{\delta}$ is the PCFB.

We next consider the case $\partial_r C>0$ and $v_0<0$.
Proposition \ref{prop:3Dl-v-} gives the PCFB.
If $\partial_r C \le \partial_r A R^{-2}$, then the condition is
\[
		\frac{(\partial_r C)^{\frac12}}{C^2}(C\partial_r A-A \partial_r C)^{\frac12}
	- \partial_r \(\frac{v_0R}{C}\)+2\partial_r \(\frac{v_0R}{C}\) \le 0.
\]
We keep the above notations $\alpha$, $\beta$, $\gamma$, and $\delta$.
Then, this is written as $\alpha \le \beta - \sqrt{\gamma}$.
Note that the right hand side is negative.
By the same argument as above, it is also written as $\alpha \le -\sqrt{\delta}$.
If $\partial_r A R^{-2}<\partial_r C \le \partial_r A (R^{-2}+v_0^2/A)$,
then the condition is
\[
		\frac{(\partial_r C)^{\frac12}}{C^2}(C\partial_r A-A \partial_r C)^{\frac12}
	 \le \left| \partial_r \(\frac{v_0R}{C}\) \right|,
\]
which is written as $\sqrt{\gamma}\le|\alpha-\beta|$.
Note that $\partial_r C < \partial_r A (R^{-2}+v_0^2/A)=C\partial_r A/A $ is equivalent to
$\gamma>0$.
By assumption, we also have $\beta<0$ and $\gamma-\beta^2=\delta-2\alpha\beta<0$.
We now show that $\alpha\ge\beta$ leads to the contradiction.
In this case, $\sqrt{\gamma}\le|\alpha-\beta|$ is equivalent to $\alpha \ge \beta + \sqrt{\gamma}$.
However, this is also written as
\begin{align*}
	0<\sqrt{\gamma}\le|\alpha-\beta| = \alpha -\beta
	& \Longleftrightarrow \alpha^2 \ge \gamma +2\alpha\beta -\beta^2=\delta>0 \\
	& \Longleftrightarrow \alpha \ge \sqrt{\delta} \text{ or } \alpha \le -\sqrt{\delta}.
\end{align*}
The last inequalities cannot be iquivalent to $\alpha \ge \beta + \sqrt{\gamma}$
since $\sqrt{\delta}>0$ and $\beta+\sqrt{\gamma}<0$.
This is the contradiction.
Hence, $\beta \ge \alpha$.
Then, $\sqrt{\gamma}\le|\alpha-\beta|=\beta-\alpha$ corresponds
to $\alpha\le -\sqrt{\delta}$.

We finally treat the case $\partial_r C \ge \partial_r A (R^{-2}+v_0^2/A)$.
We prove this condition is stronger than $\alpha \le -\sqrt{\delta}$.
An elementary computation show that
$\partial_r C \ge \partial_r A (R^{-2}+v_0^2/A)$ implies
\[
	\alpha \le \frac{C}{v_0} + \frac{v_0 R \partial_r A}{2A} < 0.
\]
Moreover, introducing the function $P(t)=\partial_r C t^2 + 2\alpha t + 2R$,
we see that
\begin{align*}
	\frac{\delta-\alpha^2}{\partial_r C} =%{}& 
	\min_{t} P(t) %\\
	\le {}& 
	P\(-\frac{v_0R}{C}\) 
	=\partial_r C \(-\frac{v_0R}{C}\)^2 + 2\alpha \(-\frac{v_0R}{C}\) + 2R \\
	={}&\frac{1}{C^2} \Bigg[ 
	\(2v_0v_0^\prime+\frac{\partial_r A}{R^2}-\frac{2A}{R^3} \)v_0^2 R^2 \\
	&{}-2(v_0+Rv_0^\prime)v_0R\(v_0^2+\frac{A}{R^2}\)+2R\(v_0^2+\frac{A}{R^2}\)^2
	\Bigg] \\
	={}& \frac{1}{C^2}\({v_0^2\partial_r A}-{2v_0v_0^\prime A}
	+\frac{2A^2}{R^3}\)\\
	={}&\frac{A}{C^2}\(\partial_r A\(R^{-2} +\frac{v_0^2}{A}\)-\partial_r C\) \le 0.
\end{align*}
\end{proof}
\subsection{PCFBs for repulsive 2D case}
We finally prove the two-dimensional case.
Though we can calculate the characteristic curve in an implicit way (\cite{ELT-IUMJ}),
we use the argument similar to the previous $n \ge 3$ case.
\begin{proof}[Proof of Proposition \ref{prop:2Dl-v+}]
We first note that $X^\prime(t,R)\ge v_0(R)>0$, $\forall t\ge0$ follows from
the same argument as in the proof of Proposition \ref{prop:3Dl-v+}.
Then, $X(t,R) \to \I$ as $t\to\I$, and, by \eqref{eq:ODE2-2}, 
\[
	\int_R^{X(t,R)} \frac{dy}{\sqrt{v_0(R)^2+A(R)\log (y/R)}} = t.
\]
for all $t\ge0$.
For simplicity, we omit the $R$ variable in the following.
Differentiate this with respect to $R$ to get
\[
	\frac{\Gamma(t)}{X^\prime(t)} - \frac{1}{v_0} -\frac12 \int_R^{X(t)} \frac{2v_0v_0^\prime-A/R + \partial_r A \log(y/R)}{(v_0^2+A\log (y/R))^{3/2}}dy = 0.
\]
We put 
\[
	B(t):= \frac{\Gamma(t)}{X^\prime(t)} =
	\frac{1}{v_0} +\frac12 \int_R^{X(t)} \frac{2v_0v_0^\prime-A/R + \partial_r A \log(y/R)}{(v_0^2+A\log (y/R))^{3/2}}dy.
\]
Since $X^\prime(t)>0$ for all $t\ge0$, $B(t)$ and $\Gamma(t)$ has the same sign.
Since $\partial_r A \ge0$ by definition, the right hand side is positive for all time if
$2v_0 v_0^\prime-A/R \ge0$.
Now, we suppose $2v_0 v_0^\prime-A/R <0$.
Recall that $X(t)\to \I$ as $t\to\I$ and that $A$ and $v_0$ are independent of time.
If $\partial_r A=0$ then one sees that there exist $t_0>0$ such that 
\[
	\int_R^{X(t_0)} \frac{ dy}{(v_0^2+A\log (y/R))^{3/2}} = \frac{2}{v_0|2 v_0 v_0^\prime-A/R|} 
\]
since $\int_R^{X(t)}(v_0^2+A\log (y/R))^{-3/2}dy\to \I$ as $t\to\I$.
This implies $\Gamma(t_0)=B(t_0)=0$, which lead to finite-time breakdown.
Let us proceed to the case $\partial_r A > 0$.
An elementary computation shows that
the minimum of $B$ is
$B\(e^{-\frac{\partial_r C}{\partial_r A}}\)$.
Therefore, under the assumption $v_0 v_0^\prime-A/R <0$ and $\partial_r A>0$,
there exists a time $t_0$ such that $\Ga(t_0)\le 0$ if and only if
\[
	\frac{1}{v_0} + \frac{1}{2} \int_R^{\exp \( -\frac{\partial_r C}{\partial_r A}\)}
	\frac{\partial_r C+\partial_r A\log y}
	{(C+A\log y)^{3/2}}dy \le 0.
\]
\end{proof}

\begin{proof}[Proof of Proposition \ref{prop:2Dl-v0}]
Let us begin with pointing out that
the exactly same argument as in the proof of Proposition \ref{prop:3Dl-v0} shows
$X^\prime(t,R)>0$ for all $t>0$ and $X(t,R)\to\I$ as $t\to\I$.
We omit $R$ variable in the following.
As in the proof of Proposition \ref{prop:3Dl-v0},
we temporarily suppose that $v_0>0$ and let $v_0\to0$ later.
Integration of \eqref{eq:ODE2-2} gives
\[
	\int_R^{X(t)} \frac{dy}{\sqrt{v_0^2+A\log (y/R)}} = t.
\]
By a change of variable $z=y/R$, the left hand side is equal to 
\[
	\int_1^{X(t)/R} \frac{R dz}{\sqrt{v_0^2+A\log z}}.
\]
Hence, differentiation with respect $R$ yields
\begin{align*}
	0=&\frac{R \partial_R(X(t)/R)}{\sqrt{v_0^2+A\log (X(t)/R)}}
	+ \int_1^{X(t)/R} \frac{dz}{\sqrt{v_0^2+A\log z}}\\
	& - R\int_1^{X(t)/R} \frac{\partial_r v_0^2+\partial_r A\log z}
	{2\(v_0^2+A\log z\)^{3/2}}dz.
\end{align*}
For a while, we omit also $t$ variable.
An elementary calculation shows
\begin{align}\label{eq:predef2Dl-v0}
	0
	=&\frac{\partial_RX}{\sqrt{v_0^2+A\log (X/R)}}
	- \frac{X}{R\sqrt{v_0^2+A\log (X/R)}}
	+ \int_1^{X/R} \frac{dz}{\sqrt{v_0^2+A\log z}} \\
	&- \frac{R\partial_r A}{2A}
	\int_1^{X/R} \frac{dz}{\sqrt{v_0^2+A\log z}}
	%\nonumber \\&
	+ \frac{Rv_0^2\partial_r A}{2A}
	\int_1^{X/R} \frac{dz}{\(v_0^2 + A\log z\)^{3/2}}
	\nonumber \\
	&- Rv_0v_0^\prime \int_1^{X/R} \frac{dz}{\(v_0^2 + A\log z\)^{3/2}}.
 	\nonumber 
\end{align}
We now show that
\begin{equation}\label{eq:limit2Dl-v0}
	\lim_{v_0 \downarrow 0} v_0 \int_1^{X/R} \frac{dz}{\(v_0^2 + A\log z\)^{3/2}} =\frac{2}{A}.
\end{equation}
Fix a small $\eps>0$.
Then, we have
\[
	\lim_{v_0 \downarrow 0} v_0 \int_{1+\eps}^{X/R}
	\frac{dz}{\(v_0^2 + A\log z\)^{3/2}} =0,
\]
since the integral is uniformly bounded with respect to $v_0$.
Moreover,
\begin{align*}
	&v_0 \int_1^{1+\eps}
	\frac{dz}{\(v_0^2 + A\log z\)^{3/2}} \\
	& \le \frac{2v_0(1+\eps)}{A}\int_1^{1+\eps}
	\frac{A}{2z(v_0^2+A\log z)^{3/2}}dz \\
	& \le \frac{2v_0(1+\eps)}{A}
	\left[\frac1{v_0}
	- \(v_0^2+A\log(1+\eps)\)^{-\frac12}\right] %\\& 
	\to \frac{2(1+\eps)}{A}
\end{align*}
as $v_0\to 0$.
Similarly, 
\begin{align*}
	v_0 \int_1^{1+\eps}
	\frac{dz}{\(v_0^2 + A\log z\)^{3/2}} 
	& \ge \frac{2v_0}{A}\int_1^{1+\eps}
	\frac{A}{2z(v_0^2+A\log z)^{3/2}}dz %\\& 
	\to \frac{2}{A}
\end{align*}
as $v_0\to 0$.
It proves \eqref{eq:limit2Dl-v0} since $\eps>0$ is arbitrary.

Taking the limit $v_0 \downarrow 0$ in \eqref{eq:predef2Dl-v0},
\begin{align*}
	0=&\frac{\partial_RX}{A^{1/2}\sqrt{\log (X/R)}}
	- \frac{X/R}{A^{1/2}\sqrt{\log (X/R)}}
	+ \frac1{A^{1/2}}\int_1^{X/R} \frac{dz}{\sqrt{\log z}} \\
	&- \frac{R\partial_r A}{2A^{3/2}}
	\int_1^{X/R} \frac{dz}{\sqrt{\log z}}
	- \frac{2Rv_0^\prime}{A} .
\end{align*}
Thus, we have
\begin{align*}
	B(t) :=\frac{\partial_R X(t)}{\sqrt{\log(X(t)/R)}} =& \frac{X(t)/R}{\sqrt{\log (X(t)/R)}}
	-  \int_1^{X(t)/R} \frac{dz}{\sqrt{\log z}}\\
	&+ \frac{R\partial_r A}{2A}\int_1^{X(t)/R} \frac{dz}{\sqrt{\log z}}
	+ \frac{2Rv_0^\prime}{A^{1/2}}.
\end{align*}
\smallbreak

{\bf Case 1.} We first assume that $\partial_r A =0$.
We put 
\[
	G(s) := \frac{s}{\sqrt{\log s}} - \int_1^{s} \frac{dz}{\sqrt{\log z}} .
\]
An elementary calculation shows
$G^\prime(s) = - (1/2)(\log s)^{-3/2} < 0$ for $s>1$,
and so $G$ is monotone decreasing.
We also see that $G^\prime$ is not integrable,
and so that $\lim_{s\to\I}G(s)=-\I$.
Since
\[
	B(t) =G(X(t)/R) + \frac{2Rv_0^\prime}{A^{1/2}},
\]
we conclude that there always exists $t_0\in(0,\I)$ such that $\Ga(t_0) = 0$.
\smallbreak

{\bf Case 2.} We next assume that $\partial_r A >0$.
We write $B(t)=:H(X(t)/R)$.
Then, it holds that
\[
	\frac{d}{ds}H(s) = - \frac{1}{2(\log s)^{3/2}}
	+\frac{R\partial_r A}{2A(\log s)^{1/2}}.
\]
Therefore, the minimum of $H$, hence of $B$, is
$H( e^{\frac{A}{R \partial_r A}})$.
The solution breaks down in finite time if and only if
this value is less than or equal to zero.
This leads to the condition
\[
	v_0^\prime R \le{}
	-\frac{\sqrt{R \partial_r A}}{2}e^{\frac{A}{R \partial_r A}}
	+\(\sqrt{A}-\frac{R\partial_r A}{2\sqrt{A}}\)  \int_0^{\sqrt{\frac{A}{R \partial_r A}}} e^{x^2}dx.
\]
\end{proof}

\begin{proof}[Proof of Proposition \ref{prop:2Dl-v-}]
If $A(R)=0$, then $X^\prime(t,R)=v_0(R)<0$ for all $t\ge 0$.
Therefore, we deduce from Lemma \ref{lem:criterion1} that the solution breaks down no latter than $t=R/|v_0(R)|$.
Hence, we assume $A(R)>0$.
Then, since $X^\prime(0,R)=v_0(R)<0$, $X^\prime(t)=-\sqrt{v_0(R)^2 + A(R) \log (X(t,R)/R)}$
as long as $X^\prime(t,R) \le 0$.
Put
\[
	t_* = \int_{Re^{-v_0^2/A}}^R \frac{dy}{\sqrt{v_0^2 + A \log (y/R)}}
	= \frac{R}{A^{1/2}e^{v_0^2/A}}\int_1^{e^{v_0^2/A}}
	\frac{dz}{\sqrt{\log z}}.
\]
Then, one sees that, for $t\in[0,t_*)$, $X(t,R)>X(t_*,R)= Re^{-v_0^2/A}>0$
and $X^\prime(t,R)<X^\prime(t_*,R)=0$.
Since $X^{\prime\prime}(t_*,R)>0$,
the same argument as in the proof of Proposition \ref{prop:3Dl-v0} shows that
$X^\prime(t,R) \ge 0$ for all $t \ge t_*$ and so that
\[
	X^\prime(t,R) =
	\begin{cases}
	- \sqrt{v_0(R)^2 + A(R) \log (X(t,R)/R)}, & \text{ for }t\le t_*, \\
	\sqrt{v_0(R)^2 + A(R) \log (X(t,R)/R)}, & \text{ for }t\ge t_*.
	\end{cases}
\]
$X(t,R)\to\I$ as $t\to\I$ is also deduced.
We omit $R$ variable in the following.
Differentiation of the identity $X(t_*,R) = Re^{-v_0^2/A}$
with respect to $R$ gives
\[
	\partial_R X(t_*) =
	e^{-v_0^2/A}\(1 - R\partial_r \(\frac{v_0^2}{A}\)\).
\]
Hence, if $R\partial_r (v_0^2/A) \ge 1$ then the solution 
breaks down no latter than $t=t_*$.
Thus, we assume $R\partial_r (v_0^2/A) < 1$ in the following.
This is equivalent to $\partial_r v_0^2 < A/R + (v_0^2/A)\partial_r A$
and to $-C/A <-\partial_r C/\partial_r A  $.
\smallbreak

{\bf Step 1.}
We first consider the condition that solution can be extended to time $t=t_*$.
For $t \le t_*$, we have
\[
	\int_{X(t)}^R \frac{dy}{\sqrt{C+A\log y}} = t.
\]
Differentiation with respect to $R$ yields
\[
	\frac{1}{\sqrt{C+A\log R}} - 
	\frac{\Ga(t)}{\sqrt{C+A\log X(t)}}
	- \frac12 \int^R_{X(t)}
	\frac{\partial_r C + \partial_r A \log y}{(C+A\log y)^{3/2}}dy = 0.
\]
For $0 \le t < t_*$, 
\[
	0 < \sqrt{C+A\log X(t)} \le \sqrt{C+ A\log R}=|v_0|
\]
holds. Therefore, $\Ga(t)$ has the same sign as
\[
	B_1(t):= \frac{\Ga(t)}{\sqrt{C+A\log X(t)}}
	=\frac{1}{|v_0|} - \frac12 \int^R_{X(t)}
	\frac{\partial_r C + \partial_r A \log y}{(C+A\log y)^{3/2}}dy.
\]
Taking time derivative, one verifies that $B_1$ takes it minimum
at $t=t_1 \in [0,t_*)$ such that 
\[
	X(t_1) = \min \(R, \exp \(-\frac{\partial_r C}{\partial_r A}\)\).
\]
Here, note that $X(t_*)=\exp(-C/A) < X(t_1)$ by assumption.
Also note that
\[
	\exp \(-\frac{\partial_r C}{\partial_r A}\) =
	R \exp \(-\frac{\partial_r v_0^2 - A/R}{\partial_r A}\).
\]
Since we have already known that $\Ga(0) =1>0$,
the solution can be extended to the time $t=t_*$ UNLESS
$\partial_r v_0^2 > A/R$ and
\[
	B_1(t_1)=
	\frac{1}{|v_0|} - \frac12 \int^R_{\exp\(-\frac{\partial_r C}{\partial_r A}\)}
	\frac{\partial_r C + \partial_r A \log y}{(C+A\log y)^{3/2}}dy
	\le 0
\]
is satisfied. Notice that this condition is a sufficient condition for finite-time breakdown.
\smallbreak

{\bf Step 2.}
We next consider the condition that the solution
can be extended from the time $t=t_*$ to $t=\I$.
For simplicity, we suppose that solutions are extended to time $t=t_*$
(we keep assuming $\partial_r v_0^2 < A/R + (v_0^2/A)\partial_r A$
holds).
Recall that, for $t > t_* $, $X^\prime(t) = \sqrt{C+A\log X(t)}>0$.

We define $t_{**}$ as a time such that $t_{**}>t_{*}$ and $X(T_{**})=R$.
Then, we have
\[
	t_{**}-t_* = \int_{Re^{-v_0^2/A}}^R \frac{dy}{\sqrt{C + A \log y}} = t_*,
\]
and so $t_{**}=2t_*$.
Thus,
\[
	\int_R^{X(t)} \frac{dy}{\sqrt{C + A \log y}} = t-2t_*
\]
for all $t \ge t_*$. As in the previous step, we set
\begin{align*}
	B_2(t):={}& \frac{\Ga(t)}{\sqrt{C+A \log X(t)}}
	=\frac{1}{|v_0|} + \frac12 \int_R^{X(t)}
	\frac{\partial_r C + \partial_r A \log y}{(C+A\log y)^{3/2}}dy
	-2\partial_r t_* \\
	={}&\frac{1}{|v_0|} + \frac12 \int_R^{X(t)}
	\frac{\partial_r v_0^2 -(A/R) + \partial_r A \log (y/R)}{(C+A\log y)^{3/2}}dy
	-2\partial_r t_*.
\end{align*}
$B_2(t)$ and $\Ga(t)$ has the same sign for $t > t_*$.
We also note that $B_2\to \I$ as $t \downarrow t_*$ because
$\Ga(t_*)>0$ and $\sqrt{C+A\log X} \to 0$ as $t \downarrow t_*$.
It holds that
\[
	\frac{d}{dt}B_2(t) = \frac{\partial_r v_0^2 -(A/R) + \partial_r A \log (X(t)/R)}{(C+A\log X(t))^{3/2}} X^\prime(t).
\]
If $\partial_r A(R)=0$ then $B_2$ is monotone decreasing
by assumption $\partial_r v_0^2 - A/R <0$.
Moreover, $\frac{d}{dt}B_2(t)$ is uniformly bounded by $(\partial_r v_0^2 -(A/R))/|v_0|<0$
from above, and so there exists time $t_2$ such that $B_2(t_2) = 0$.
Therefore, now we suppose $\partial_r A(R) >0$.

$B_2$ takes it minimum at
$t=t_2$ such that $X(t_2)=\exp(-{\partial_r C}/{\partial_r A})$.
Therefore, the solution can be extended to $t=\I$ if and only if
\[
	B_2(t_2)=\frac{1}{|v_0|}+ \frac12 \int_R^{\exp\(-\frac{\partial_r C}{\partial_r A}\)}
	\frac{\partial_r C + \partial_r A \log y}{(C+A\log y)^{3/2}}dy
	-2\partial_r t_* >0.
\]
\end{proof}

\subsection*{Acknowledgments}
The author would like to thank the referee for 
reading the manuscript very carefully and
giving many valuable suggestions.
The author also expresses his deep gratitude
to Professor Yoshio Tsutsumi for his 
valuable advice and constant encouragement.
This research is supported by JSPS fellow.

\bibliographystyle{amsplain}
\bibliography{caustic}
\end{document}